\documentclass[graybox]{svmult}

\usepackage{pict2e}

\usepackage{mathptmx}       
\usepackage{helvet}         
\usepackage{courier}        
\usepackage{type1cm}        
\usepackage{makeidx}         
\usepackage{graphicx}        
\usepackage{multicol}        
\usepackage[bottom]{footmisc}

\RequirePackage[OT1]{fontenc}
\RequirePackage{amsmath,amssymb}

\numberwithin{equation}{section}

\newtheorem{algorithm}{Algorithm}[section]

\newcommand{\indicator}[2]{\delta\left(#1 \mid #2 \right)}

\newcommand{\B}[1]{{\bf #1}}
\newcommand{\Sc}[1]{{\mathcal{#1}}}
\newcommand{\R}[1]{{\rm #1}}
\newcommand{\mB}[1]{{\mathbb{#1}}}
\newcommand{\eR}{\overline{\mathbb{R}}}

\newcommand{\Del}{\Delta}
\newcommand{\Rn}{\mB{R}^n}

\makeindex             

\begin{document}

\title*{Optimization viewpoint on Kalman smoothing, with applications to robust and sparse estimation. }
\titlerunning{Optimization viewpoint on Kalman smoothing}
\author{Aleksandr Y. Aravkin and James V. Burke and Gianluigi Pillonetto}
\institute{Aleksandr Y. Aravkin \at IBM T.J. Watson Research Center\\Yorktown Heights, NY, 10598\\ \email{saravkin@us.ibm.com}
\and James V. Burke \at Department of Mathematics,\\ University of Washington, Seattle, WA\\ \email{jvburke@uw.edu}
\and Gianluigi Pillonetto \at Control and Dynamic Systems Department of Information Engineering, \\University of Padova, Padova, Italy,\\ \email{giapi@dei.unipd.it}
}
%
%
\maketitle

\abstract*{In this chapter, we present the optimization formulation of the Kalman filtering and smoothing problems, and use this perspective 
to develop a variety of extensions and applications. We first formulate 
classic Kalman smoothing as a least squares problem, highlight special structure, and show that the classic filtering and smoothing 
algorithms are equivalent to a particular algorithm for solving this problem. Once this equivalence is established, 
we present extensions of Kalman smoothing to systems with nonlinear process and measurement models, 
systems with linear and nonlinear inequality constraints, systems with outliers in the measurements or sudden changes in the state, 
and systems where the sparsity of the state sequence must be accounted for. All extensions preserve the computational 
efficiency of the classic algorithms, and most of the extensions are illustrated with numerical 
examples, which are part of an open source Kalman smoothing Matlab/Octave package.}

\abstract{In this chapter, we present the optimization formulation of the Kalman filtering and smoothing problems, and use this perspective 
to develop a variety of extensions and applications. We first formulate 
classic Kalman smoothing as a least squares problem, highlight special structure, and show that the classic filtering and smoothing 
algorithms are equivalent to a particular algorithm for solving this problem. Once this equivalence is established, 
we present extensions of Kalman smoothing to systems with nonlinear process and measurement models, 
systems with linear and nonlinear inequality constraints, systems with outliers in the measurements or sudden changes in the state, 
and systems where the sparsity of the state sequence must be accounted for. All extensions preserve the computational 
efficiency of the classic algorithms, and most of the extensions are illustrated with numerical 
examples, which are part of an open source Kalman smoothing Matlab/Octave package.}

\section{Introduction}

Kalman filtering and smoothing methods form a broad category 
of computational algorithms used for inference on 
noisy dynamical systems. Over the last fifty years,
these algorithms have become a gold standard in 
a range of applications, including space exploration, 
missile guidance systems, general 
tracking and navigation, and weather prediction. 
In 2009, Rudolf Kalman received the 
National Medal of Science from President Obama for the invention
of the Kalman filter. Numerous books and papers 
have been written on these methods and their extensions, 
addressing modifications for use in nonlinear systems, 
smoothing data over time intervals, 
improving algorithm robustness to bad measurements, 
and many other topics.

\begin{figure}
\begin{center}
{\includegraphics[width=5in]{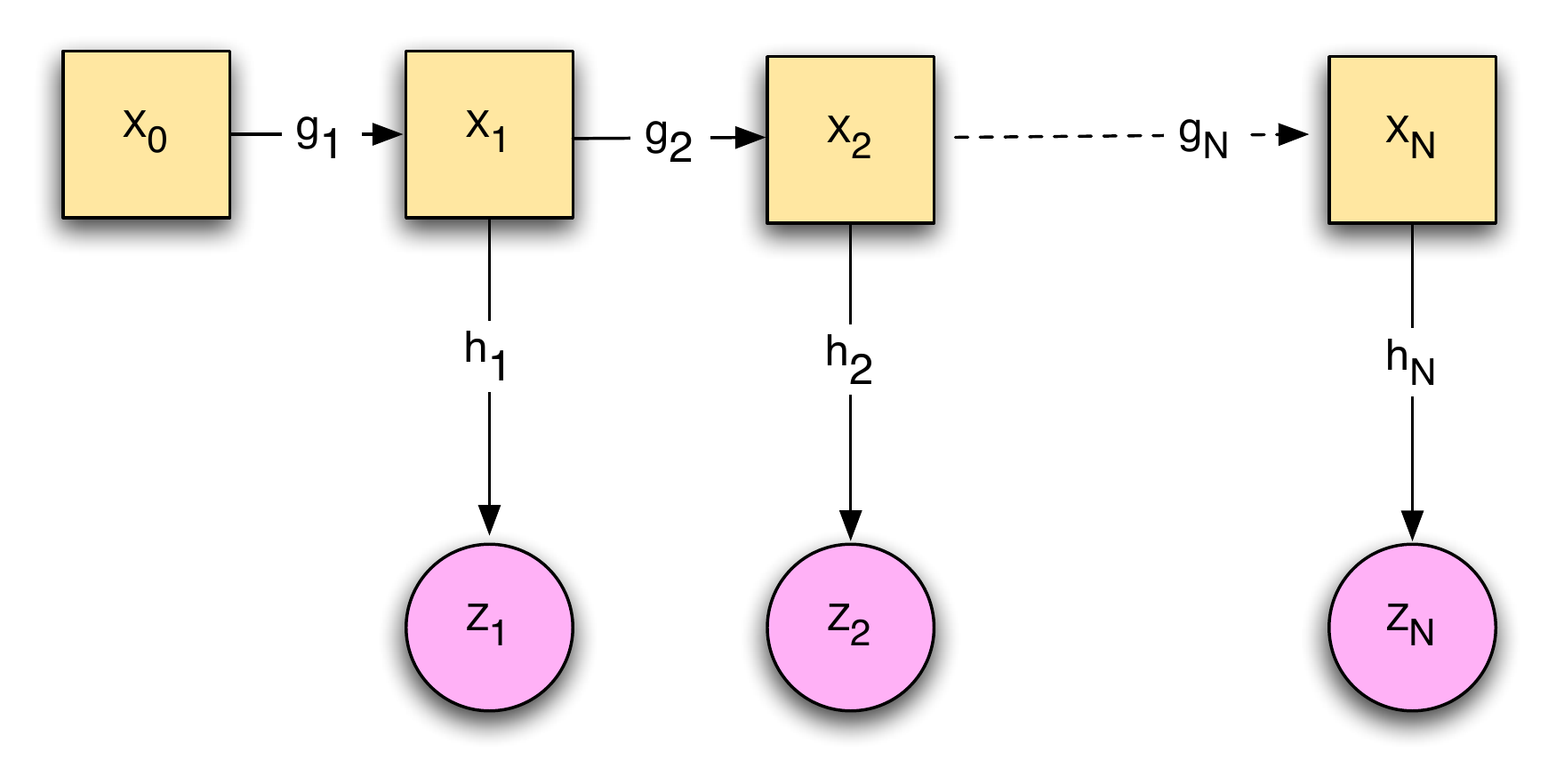}}
\end{center}
\caption{
	\label{kalmanFramework}
Dynamic systems amenable to Kalman smoothing methods. 
}
\end{figure}

The classic Kalman filter~\cite{kalman} is almost always presented as a set of recursive equations, and the
classic Rauch-Tung-Striebel (RTS) fixed-interval smoother~\cite{RTS} is typically formulated as two coupled Kalman
filters. An elegant derivation based on projections onto spaces spanned by random variables can be found in \cite{Ansley}.
In this chapter, we use the terms `Kalman filter' and `Kalman smoother' 
much more broadly, {\bf including any method 
of inference on any dynamical system fitting the graphical representation of Figure~\ref{kalmanFramework}.}
Specific mathematical extensions we consider include 
\begin{itemize}

\item Nonlinear process and measurement models. 
\item Inequality state space constraints.
\item Different statistical models for process and measurement errors. 
\item Sparsity constraints. 
\end{itemize}
We also show numerous applications of these extensions. 

The key to designing tractable inference methods for the above applications 
is an optimization viewpoint, which we develop in the classic Kalman smoothing case
and then use to formulate and solve {\it all} of the above extensions. 
Though it has been known for many years that the Kalman filter provides the 
 maximum {\it a posteriori}  estimate for linear systems subject to Gaussian noise,  
 the optimization perspective underlying this idea has not been fully deployed 
 across engineering applications. 
Notably, several groups (starting in 1977) 
have discovered and used variants of this perspective 
to implement extensions to Kalman filtering 
and smoothing, including singular filtering 
(\cite{Paige1977, Paige1981, Paige1985}),
robust smoothing (\cite{Fahr1998, Aravkin2011tac}), 
nonlinear smoothing with inequality 
state space constraints (\cite{Bell1994, Bell2008}), and 
sparse Kalman smoothing~\cite{Angelosante2009}.

We focus exclusively on smoothing here, leaving online 
applications of these ideas to future work (see~\cite{Pillonetto2010} for an 
example of using a smoother for an online application). 
We start by presenting the classic RTS smoothing algorithm in Section~\ref{Classic},
and show that the well-known recursive equations 
are really an algorithm to solve a least squares system with special structure. 
Once this is clear, it becomes much easier to discuss novel extensions, since 
as long as special structure is preserved, their computational
cost is on par with the classic smoother (or, put another way, 
the classic smoothing equations are viewed as a particular way to solve
key subproblems in the extended approaches). 

In the subsequent sections, we build novel extensions, 
briefly review theory, discuss the special structure, and 
present numerical examples for a variety of applications. In Section~\ref{Nonlinear}, we 
formulate the problem for smoothing with nonlinear process 
and measurement models, and show how to solve it. 
In Section~\ref{Constrained}, we show how state space 
constraints can be incorporated, and the resulting problem solved
using interior point techniques. 
In Section~\ref{Robust}, we review two recent Kalman smoothing
formulations that are highly robust to measurement errors. 
Finally, in Section~\ref{Sparse}, we review recent work in 
sparse Kalman smoothing, and show how sparsity can be 
incorporated into the other extensions. 
We end the chapter with discussion in Section~\ref{Conclusions}.

\section{Optimization Formulation and RTS Smoother}
\label{Classic}

\subsection{Probabilistic model}
The model corresponding to 
Figure \ref{kalmanFramework} is specified as follows: 

%
\begin{equation}
\label{IntroGaussModel}
\begin{array}{rcll}
	\B{x_1}&= & g_1(x_0)+\B{w_1},
	\\
	\B{x_k} & = & g_k (\B{x_{k-1}})  + \B{w_k}& k = 2 , \ldots , N,
	\\
	\B{z_k} & = & h_k (\B{x_k})      + \B{v_k}& k = 1 , \ldots , N\;,
\end{array}
\end{equation}
where $\B{w_k}$, $\B{v_k}$ are mutually independent random variables
with known positive definite covariance matrices $Q_k$ and $R_k$, respectively. 
We have $\B{x_k}, \B{w_k} \in \mB{R}^n$, and $\B{z_k}, \B{v_k} \in \mB{R}^{m(k)}$\;, so 
measurement dimensions can vary between time points. 
The classic case is obtained by making the following assumptions:
\begin{enumerate}
\item 
\label{linearAss}
$x_0$ is known, and $g_k$, $h_k$ are known {\it linear} functions, which we denote by 
\begin{equation}\label{IntroLinearModel}
\begin{aligned}
g_k(x_{k-1}) &= G_kx_{k-1} 
\quad h_k(x_k) &= H_k x_k 
\end{aligned}
\end{equation}
where $G_k \in \mB{R}^{n\times n}$ and $H_k \in \mB{R}^{m(k)\times n}$, \item 
\label{gaussianAss}
$\B{w_k}$, $\B{v_k}$ are mutually independent {\it Gaussian} random variables.
\end{enumerate}
In later sections, we will show how to relax these classic assumptions, 
and what gains can be achieved once they are relaxed. 
In this section, we will formulate estimation of the {\it entire} state sequence,
$x_1,x_2,\dots,x_N,$
as an optimization problem, and show how the RTS smoother solves it. 

\subsection{Maximum {\it a posteriori}  formulation}
\label{sec:MAP}
To begin, we formulate the maximum {\it a posteriori} (MAP) problem under 
linear and Gaussian assumptions. 
Using Bayes' theorem, we have 
\begin{equation}\label{MAP}
\begin{aligned}
P\left(\{x_k\} \big| \{z_k\}\right) &\propto P\left(\{z_k\}\big|\{x_k\}\right)P\left(\{x_k\}\right)\\
& = \prod_{k=1}^N P(\{v_k\})  P(\{w_k\})\\
& \propto 
\prod_{k=1}^N \exp\Big( -\frac{1}{2}(z_k - H_k x_k)^\top R_k^{-1}(z_k - H_k x_k) \\
& \quad\quad-\frac{1}{2}(x_k - G_k x_{k-1})^\top Q_k^{-1}(x_k - G_k x_{k-1})\Big)\;.
\end{aligned}
\end{equation}
A better (equivalent) formulation to~\eqref{MAP} is minimizing its negative log posterior: 
\begin{equation}
\label{logMAP}
\begin{aligned}
\min_{\{x_k\}} f(\{x_k\})&:= 
\sum_{k=1}^N\frac{1}{2}(z_k - H_k x_k)^\top R_k^{-1}(z_k - H_k x_k) +\frac{1}{2}(x_k - G_k x_{k-1})^\top Q_k^{-1}(x_k - G_k x_{k-1})\;.
\end{aligned}
\end{equation}

To simplify the problem, we now introduce data structures that capture the entire state sequence, measurement sequence, 
covariance matrices, and initial conditions. 

Given a sequence of column vectors $\{ u_k \}$
and matrices $ \{ T_k \}$ we use the notation
\[
\R{vec} ( \{ u_k \} )
=
\begin{bmatrix}
u_1 \\ u_2  \\ \vdots \\ u_N
\end{bmatrix}
\; , \;
\R{diag} ( \{ T_k \} )
=
\begin{bmatrix}
T_1    & 0      & \cdots & 0 \\
0      & T_2    & \ddots & \vdots \\
\vdots & \ddots & \ddots & 0 \\
0      & \cdots & 0      & T_N
\end{bmatrix} .
\]
We now make the following definitions:
\begin{equation}\label{defs}
\begin{aligned}
R       & =  \R{diag} ( \{ R_k \} )
\\
Q       & =  \R{diag} ( \{ Q_k \} )
\\
H       & = \R{diag} (\{H_k\} )
\end{aligned}\quad \quad
\begin{aligned}
x       & = \R{vec} ( \{ x_k \} )
\\
w      &  = \R{vec} (\{g_0, 0, \dots, 0\})
\\
z      & = \R{vec} (\{z_1,  z_2, \dots, z_N\})
\end{aligned} \quad \quad
\begin{aligned}
G  & = \begin{bmatrix}
    \R{I}  & 0      &          &
    \\
    -G_2   & \R{I}  & \ddots   &
    \\
        & \ddots &  \ddots  & 0
    \\
        &        &   -G_N  & \R{I}
\end{bmatrix}\;,
\end{aligned}
\end{equation}
where $g_0:=g_1(x_0)=G_1x_0$.

With definitions in~\eqref{defs}, problem~\eqref{logMAP} can be written
\begin{equation}\label{fullLS}
\min_{x} f(x) =  \frac{1}{2}\|Hx - z\|_{R^{-1}}^2 + \frac{1}{2}\|Gx -w\|_{Q^{-1}}^2\;,
\end{equation}
where $\|a\|_M^2 = a^\top Ma$.
We knew the MAP  was a least squares problem already, but now the structure is fully transparent. 
In fact, we can write down the closed form solution by taking the gradient
of~\eqref{fullLS} and setting it equal to $0$:
\[
\begin{aligned}
0 & = H^\top R^{-1}(Hx-z) + G^\top Q^{-1}(Gx - w) 
\\
& =  (H^\top R^{-1} H + G^\top Q^{-1} G) x  - H^\top R^{-1}z -G^\top Q^{-1}w\;.
\end{aligned}
\]
The smoothing estimate is therefore given by solving the linear system
\begin{equation}\label{smoothingSol}
(H^\top R^{-1} H + G^\top Q^{-1} G) x =  H^\top R^{-1}z + G^\top Q^{-1}w\;.
\end{equation}

\subsection{Special subproblem structure}

The linear system in~\eqref{smoothingSol} has a very special structure: it is 
a symmetric positive definite block tridiagonal matrix. This can be immediately observed 
from the fact that both $G$ and $Q$ are positive definite. 
 To be specific, it is given by
\begin{equation}
\label{hessianApprox}
C = (H^\top R^{-1} H + G^\top Q^{-1} G)
=
\begin{bmatrix}
C_1 & A_2^\R{T} & 0 & \\
A_2 & C_2 & A_3^\R{T} & 0 \\
0 & \ddots & \ddots& \ddots & \\
& 0 & A_N & C_N
\end{bmatrix} ,
\end{equation}
with $A_k \in \mB{R}^{n\times n}$ and
$C_k \in \mB{R}^{n\times n}$ defined as follows:
\begin{eqnarray}
\label{KalmanData}
\nonumber
A_k
&=&
-Q_k^{-1}G_{k}\; , \;\\
\nonumber
C_k
&=&
Q_k^{-1} + G_{k+1}^\top Q^{-1}_{k+1}G_{k+1} +H_k^\top R_k^{-1}H_k\; .\\
\end{eqnarray}

The special structure of the matrix $C$ in~\eqref{hessianApprox} can be exploited
to solve the linear system equivalent to the Kalman smoother. While a structure-agnostic   
matrix inversion scheme has complexity $O(n^3N^3)$, 
exploiting the block tridiagonal structure reduces this complexity to $O(n^3N)$.  

A straightforward algorithm for solving any symmetric positive definite block tridiagonal linear system 
is given in~\cite{Bell2000}. We review it here, since it is essential to build the connection
to the standard viewpoint of the RTS smoother. 

\subsection{Block tridiagonal (BT) algorithm}

Suppose for \( k = 1 , \ldots , N \),
\( c_k \in \B{R}^{n \times n} \),
\( e_k \in \B{R}^{n \times \ell} \),
\( r_k \in \B{B}^{n \times \ell} \),
and for \( k = 2 , \ldots , N \),
\( a_k \in \B{R}^{n \times n} \). 
We define the corresponding block tridiagonal system of equations
\begin{equation}
\label{BlockTridiagonalEquation}
\small
\left( \begin{matrix}
c_1     & a_2^\R{T}  & 0       & \cdots        & 0      \\
a_2     & c_2        &         &               & \vdots \\
\vdots  &            & \ddots  &               & 0         \\
0       &            & a_{N-1} & c_{N-1}       & a_N^\R{T} \\
0       & \cdots     & 0       & a_N           & c_N 
\end{matrix} \right)
\left( \begin{array}{c} 
	e_1 \\ e_2 \\ \vdots \\ e_{N-1} \\ e_N
\end{array} \right)
=
\left( \begin{matrix} 
	r_1 \\ r_2 \\ \vdots \\  r_{N-1} \\ r_N
\end{matrix} \right)
\end{equation}

The following algorithm for \eqref{BlockTridiagonalEquation} is given in \cite[Algorithm 4]{Bell2000}. 

\begin{algorithm}
\label{blockAlg}
The inputs to this algorithm are 
\( \{ a_k \} \),
\( \{ c_k \} \),
and
\( \{ r_k \} \).
The output is a sequence \( \{ e_k \} \)
that solves equation~(\ref{BlockTridiagonalEquation}). 
\end{algorithm}
\begin{enumerate}

\item
Set \( d_1 = c_1  \) and \( s_1 = r_1 \).

\item
For \( k = 2, \ldots , N \),
set \( d_k = c_k - a_{k}^\R{T} d_{k-1}^{-1} a_{k} \), ~
\( s_k = r_k - a_{k}^\R{T} d_{k-1}^{-1} s_{k-1} \).

\item
Set \( e_N = d_N^{-1} s_N \).

\item
For \( k = N-1 , \ldots , 1 \),
set \( e_k = d_k^{-1} ( s_k - a_{k+1} e_{k+1} ) \).


\end{enumerate}
Note that after the first two steps of Algorithm~\ref{blockAlg}, we have arrived at a linear
system equivalent to~\eqref{BlockTridiagonalEquation} but upper triangular: 
\begin{equation}
\label{BlockTridiagonalMod}
\small
\left( \begin{matrix}
d_1     & a_2^\R{T}  & 0       & \cdots        & 0      \\
0     & d_2        &         &               & \vdots \\
\vdots  &            & \ddots  &               & 0         \\
0       &            & 0 & d_{N-1}       & a_N^\R{T} \\
0       & \cdots     & 0       & 0           & d_N 
\end{matrix} \right)
\left( \begin{array}{c} 
	e_1 \\ e_2 \\ \vdots \\ e_{N-1} \\ e_N
\end{array} \right)
=
\left( \begin{matrix} 
	s_1 \\ s_2 \\ \vdots \\  s_{N-1} \\ s_N
\end{matrix} \right)
\end{equation}
The last two steps of the algorithm then simply back-solve for the $e_k$. 

\subsection{Equivalence of Algorithm~\eqref{blockAlg} to Kalman Filter and RTS Smoother}
\label{equivalence}
Looking at the very first block, we now substitute in the Kalman data structures~\eqref{KalmanData}
into step 2 of Algorithm~\ref{blockAlg}: 
\begin{equation}\label{Equivalence}
\small
\begin{aligned}
d_2 &= c_2 - a_{2}^\R{T} d_{1}^{-1} a_{2}\\
& = \underbrace{\underbrace{Q_2^{-1}   - \left(Q_2^{-1}G_{2}\right)^\top\left(\underbrace{Q_1^{-1} + H_1^\top R_1^{-1}H_1}_{P_{1|1}^{-1}} + G_{2}^\top Q^{-1}_{2}G_{2} \right )^{-1} \left(Q_2^{-1}G_{2}\right)}_{P_{2|1}^{-1}} + H_2^\top R_2^{-1}H_2}_{P_{2|2}^{-1}}
+ G_{3}^\top Q^{-1}_{3}G_{3}
\end{aligned}
\end{equation}
These relationships can be seen quickly from~\cite[Theorem 2.2.7]{AravkinThesis2010}. 
The matrices $P_{k|k}$, $P_{k|k-1}$ are common to the Kalman filter framework: they represent 
covariances of the state at time $k$ given the the measurements $\{z_1, \dots, z_k\}$, and the covariance of the 
a priori state estimate at time $k$ given measurements $\{z_1, \dots, z_{k-1}\}$, respectively.

From the above computation, we see that 
\[
d_2 = P_{2|2}^{-1} +G_{3}^\top Q^{-1}_{3}G_{3}\;. 
\]
By induction, it is easy to see that in fact 
\[
d_k = P_{k|k}^{-1} +G_{k+1}^\top Q^{-1}_{k+1}G_{k+1}\;. 
\]
We can play the same game with $s_k$. Keeping in mind that  
$r = H^\top R^{-1}z + G^\top Q^{-1}w$, we have
\begin{equation}\label{EquivalenceRHS}
\small
\begin{aligned}
s_2 &= r_2 - a_{2}^\R{T} d_{1}^{-1} r_{1}\\
& = \underbrace{
H_2^\top R_2^{-1}z_2 + \underbrace{\left(Q_2^{-1}G_{2}\right)^\top\left(\underbrace{Q_1^{-1} + H_1^\top R_1^{-1}H_1}_{P_{1|1}^{-1}} + G_{2}^\top Q^{-1}_{2}G_{2} \right )^{-1}
\left(H_1^\top R_1^{-1}z_1 + G_1^\top P_{0|0}^{-1}x_0 \right)}_{a_{2|1}}}_{a_{2|2}}
\end{aligned}
\end{equation}
These relationships also follow from~\cite[Theorem 2.2.7]{AravkinThesis2010}. 
The quantities $a_{2|1}$ and $a_{2|2}$ are from the information filtering literature, and are less commonly known: they are preconditioned estimates 
\begin{equation}
\label{infoStructures}
\begin{aligned}
a_{k|k} &= P_{k|k}^{-1}x_{k}\;, \quad a_{k|k-1} & = P_{k|k-1}^{-1}x_{k|k-1}\;.
\end{aligned}
\end{equation}
Again, by induction we have precisely that $s_k = a_{k|k}$. 

When you put all of this together, you see that step 3 of Algorithm~\ref{blockAlg} is given by 
\begin{equation}\label{KalmanFilter}
\begin{aligned}
e_N &= d_N^{-1}s_N  
 = \left(P_{N|N}^{-1} + 0 \right)^{-1} P_{N|N}^{-1}x_{k|k}
 = x_{k|k}\;,
\end{aligned}
\end{equation}
so in fact $e_N$ is the Kalman filter estimate (and the RTS smoother estimate) for time point $N$.  

Step $4$ of Algorithm~\ref{blockAlg} then implements the backward Kalman filter, computing the 
smoothed estimates $x_{k|N}$ by back-substitution. {\bf Therefore the RTS smoother is 
Algorithm~\ref{blockAlg} applied to~\eqref{smoothingSol}}. 

The consequences are profound --- instead of working with the kinds expressions seen in~\eqref{EquivalenceRHS}
and~\eqref{Equivalence}, we can think at a high level, focusing on~\eqref{fullLS}, and simply 
using Algorithm~\ref{blockAlg} (or variants) as a subroutine.  As will become apparent, the key to all extensions 
is preserving the block tridiagonal structure in the subproblems, so that Algorithm~\ref{blockAlg} can be used. 

\subsection{Numerical Example: Tracking a Smooth Signal}
\label{sec:LinearExample}

\begin{figure}\label{fig:linearSmooth}
\begin{center}
{\includegraphics[width=4in]{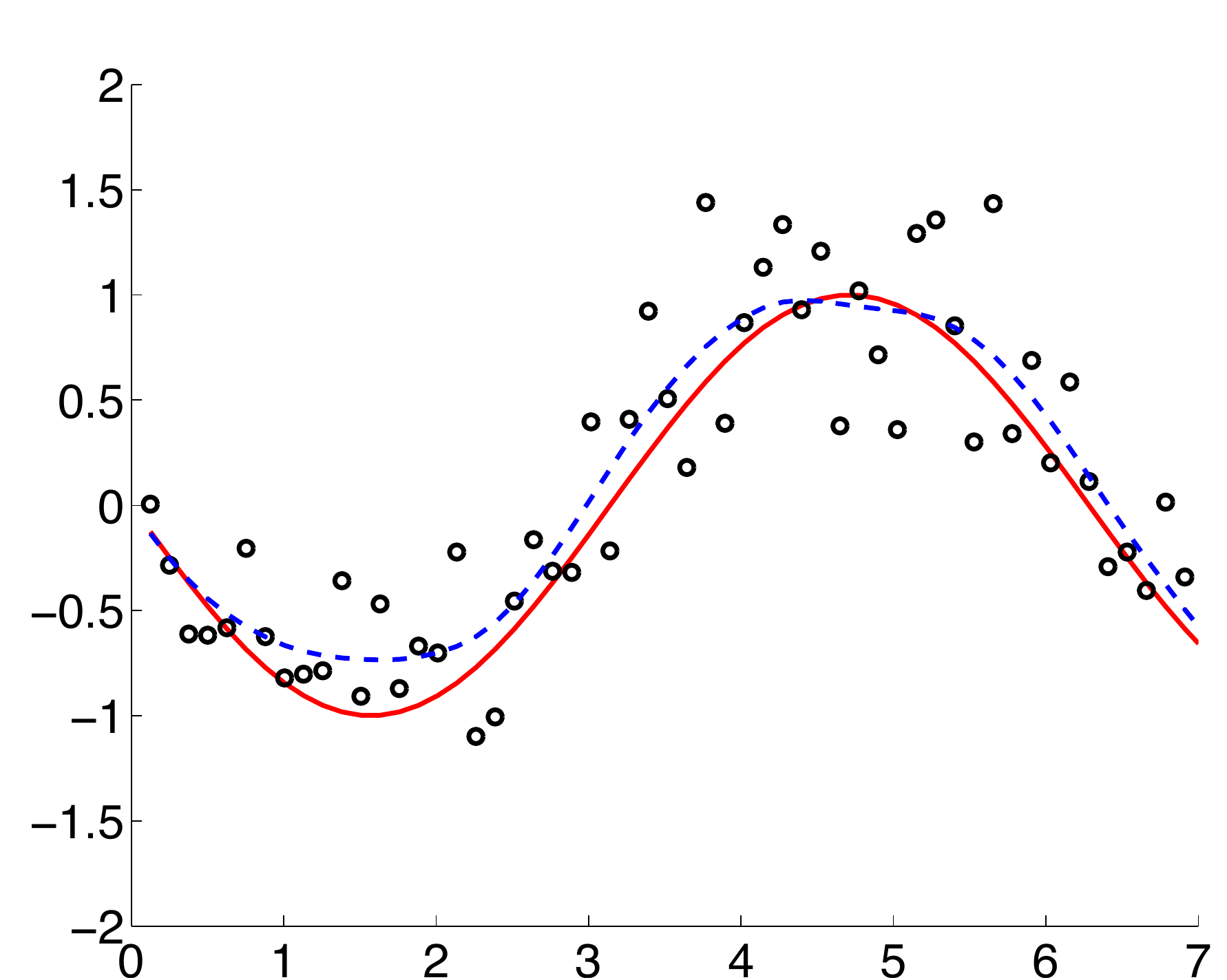}}
\end{center}
\caption{
Tracking a smooth signal (sine wave) using a generic linear process model~\eqref{linearSmooth} and direct (noisy) measurements~\eqref{directMeasurement}. 
Red solid line is true signal, blue dashed line is Kalman (RTS) smoother estimate. Measurements are displayed as circles. 
}
\end{figure}

In this example, we focus on a very useful and simple model: the process model for a {\it smooth} signal. 
Smooth signals arise in a range of applications: physics-based models, biological 
data, and financial data all have some inherent smoothness. 

A surprisingly versatile technique for modeling {\it any} such process is to 
treat it as integrated Brownian motion. We illustrate on a scalar time series $x$. 
We introduce a new derivative state $\dot x$, with process model 
\(
\dot x_{k+1} = \dot x_{k} + \dot w_k\;,
\)
and then model the signal $x$ or interest as 
\(
x_{k+1} = x_k + \dot x_{k} \Delta t + w_k\;.
\)
Thus we obtain an augmented (2D) state with process model 
\begin{equation}
\label{linearSmooth}
\begin{bmatrix}
\dot x_{k+1} \\ x_{k+1}
\end{bmatrix} 
= 
\begin{bmatrix}
I & 0 
\\
\Delta t & I
\end{bmatrix}
\begin{bmatrix}
\dot x_{k} \\ x_{k}
\end{bmatrix} 
+
\begin{bmatrix}
\dot w_{k} \\ w_{k}
\end{bmatrix} \;.
\end{equation}

Using a well-known connection to stochastic differential equations (see \cite{Jaz,Oks,Bell2008}) 
we use covariance matrix 
\begin{equation}\label{SDEvar}
Q_k = \sigma^2 \begin{bmatrix}\Delta t & \Delta t^2/2 \\ \Delta t^2 / 2 & \Delta t^3/ 3 \end{bmatrix}\;.
\end{equation}

Model equations~\eqref{linearSmooth} and~\eqref{SDEvar} can be applied as a process model for any 
smooth process.  For our numerical example, we take direct measurements of the $\sin$ function, which is very smooth. 
Our measurement model therefore is 
\begin{equation}\label{directMeasurement}
z_k = H_k x_k + v_k\;, \quad H_k  = \begin{bmatrix} 0 & 1 \end{bmatrix}\;.
\end{equation}

The resulting fit is shown in Figure~\ref{fig:linearSmooth}. The measurements guide the estimate to the true smooth time series, 
giving very nice results. The figure was generated 
using the \verb{ckbs{ package~\cite{ckbs}, specifically using the example file \verb{affine_ok.m{.
Measurement errors were generated using $R_k = .35^2$, and this value was given to the smoother. 
The $\sigma^2$ in~\eqref{SDEvar} was taken to be $1$.
The program and example are available for download from COIN-OR.

\section{Nonlinear Process and Measurement Models}
\label{Nonlinear}

In the previous section, we have shown that when $g_k$ and $h_k$
in model~\eqref{IntroGaussModel} are linear, and $\B{v_k, w_k}$
are Gaussian, then the smoother is equivalent to solving a least squares
problem~\eqref{fullLS}. We have also shown 
that the filter estimates appear as intermediate results when 
one uses Algorithm~\ref{blockAlg} to solve the problem. 

In this section, we turn to the case where $g_k$ and $h_k$ are nonlinear. 
We first formulate the smoothing problem as a maximum {\it a posteriori} (MAP) 
 problem, and show that it is a nonlinear least squares (NLLS)
problem. To set up for later sections, we also introduce the broader class of {\it convex composite} problems. 

We then review the standard Gauss-Newton method in the broader context 
of convex composite models, and show that when applied to the NLLS problem,  
each iteration is equivalent to solving~\eqref{fullLS}, and therefore 
to a full execution of the RTS smoother. 
We also show how to use a simple line search to guarantee convergence 
of the method to a local optimum of the MAP problem. 

This powerful approach, known for at least 20 years~\cite{Fahr1991,Bell1993,Bell1994}, 
is rarely used in practice; instead practitioners favor the EKF or the UKF~\cite{merwe,JUW}, 
neither of which converge to a (local) MAP solution. 
MAP approaches work very well for a broad range of applications, 
and it is not clear why one would throw away an efficient MAP solver in favor of another scheme. 
To our knowledge, the optimization (MAP) approach has never been 
included in a performance comparison of `cutting edge' methods, such as~\cite{Lef}. 
While such a comparison is not in the scope of this work, we lay the foundation 
by providing a straightforward exposition of the optimization approach and a reproducible numerical illustration (with publicly available code)
for smoothing the Van Der Pol oscillator, a well known problem where the process model is a nonlinear ODE.

\subsection{Nonlinear Smoother Formulation and Structure}
\label{sec:ConvexComposite}

In order to develop a notation analogous to~\eqref{fullLS}, we define
functions $g:\mathbb{R}^{nN}\rightarrow\mathbb{R}^{n(N+1)}$ and 
$h:\mathbb{R}^{nN} \rightarrow \mathbb{R}^{M}$, with $M = \sum_k m_k$,
from components $g_k$ and $h_k$ as follows. 
\begin{equation}
\label{ghDef}
g(x) = \begin{bmatrix}x_1 \\ x_2 - g_2(x_1) \\ \vdots \\ x_N - g_N(x_{N-1}) \end{bmatrix}\;,
\quad
h(x) = \begin{bmatrix} h_1(x_1) \\ h_2(x_2) \\ \vdots \\ h_N(x_N)\end{bmatrix}\;.
\end{equation}

With this notation, the MAP problem, obtained exactly as in Section~\ref{sec:MAP}, is given by 
\begin{equation}
\label{fullNLLS}
\min_x f(x) = \frac{1}{2}\|g(x) -w\|_{Q^{-1}}^2 + \frac{1}{2}\|h(x) - z\|_{R^{-1}}^2\;,
\end{equation}
where $z$ and $w$ are exactly as in~\eqref{defs}, so that $z$ is the {\it entire} vector 
of measurements, and $w$ contains the initial estimate $g_1(x_0)$ in the first $n$ entries, 
and zeros in the remaining $n(N-1)$entries. 

We have formulated the nonlinear smoothing problem as a 
nonlinear least-squares (NLLS) problem --- compare~\eqref{fullNLLS} with~\eqref{fullLS}. 
We take this opportunity to note that NLLS problems are a special example of a more general structure. 
Objective~\eqref{fullNLLS} may be written as a composition of a convex function $\rho$ with a smooth function $F$:
\begin{equation}\label{convexComposite}
f(x) = \rho(F(x))\;, 
\end{equation}
where
\begin{equation}
\label{CCparts}
\rho\left(\begin{matrix} y_1\\ y_2\end{matrix}
\right) 
= \frac{1}{2}\|y_1\|_{Q^{-1}}^2 + \frac{1}{2}\|y_2\|_{R^{-1}}^2\;, \quad F(x) = \begin{bmatrix} g(x) -w \\ h(x) - z\end{bmatrix}\;.
\end{equation}

As we show in the next sub-section, problems of general form~\eqref{convexComposite} can be solved using the Gauss-Newton method,
which is typically associated specifically with NLLS problems. Presenting the Gauss-Newton right away in the more general setting will make it 
easier to understand extensions in the following sections of the chapter. 

\subsection{Gauss-Newton Method for Convex Composite Models}
\label{sec:GaussNewton}

The Gauss-Newton method can be used to solve problems of the form~\eqref{convexComposite},
and it uses a very simple strategy: {\bf iteratively linearizing the smooth function $F$}
\cite{Burke85}. 
More specifically, the Gauss-Newton method is an iterative method of the form
\begin{equation}\label{GN}
x^{\nu + 1} = x^\nu + \gamma^\nu d^\nu\;,
\end{equation}
where $d^\nu$ is the Gauss-Newton search direction, and $\gamma^\nu$ is a scalar that guarantees
\begin{equation}\label{descent}
f(x^{\nu+1}) < f(x^\nu)\;.
\end{equation}
The direction $d^\nu$ is obtained by solving the subproblem 
\begin{equation}\label{GNsub}
d^\nu = \arg\min_{d} \tilde f(d) := \rho\left(F(x^\nu) + \nabla F(x^\nu)^\top d\right)\;.
\end{equation}
We then set 
\[
\tilde\Del f(x^\nu)= \tilde f(d^\nu)-f(x^\nu).
\]

By \cite[Lemma 2.3, Theorem 3.6]{Burke85}, 
\begin{equation}\label{dirDer}
f'(x^\nu; d^\nu) \leq \tilde\Del f(x^\nu) \le 0\;,
\end{equation}
with equality if and only if $x^\nu$ is a first-order stationary point for $f$.
This implies that a suitable stopping criteria for the algorithm is the
condition $\Del f(x^\nu)\sim 0$. Moreover, 
$x^\nu$ is not a first-order stationary point for $f$, then 
the direction $d^\nu$ is a direction of strict descent for $f$ at $x^\nu$. 

Once the direction $d^\nu$ is obtained with $\tilde\Del f(x^\nu) < 0$, a step-size $\gamma^\nu$ 
is obtained by a standard backtracking
line-search procedure: pick  
a values $0<\lambda < 1$ and $0<\kappa<1$ (e.g., $\lambda = 0.5$ and $\kappa=0.001$) 
and evaluate 
$f(x^\nu +  \lambda^s d^\nu)$, $s = 0, 1, 2, \dots,$ until
\begin{equation}\label{armijo ineq}
f(x^\nu +  \lambda^s d^\nu)\le f(x^\nu)+\kappa\lambda^s\tilde\Del f(x^\nu)
\end{equation}
is satisfied for some 
$\bar s$, 
then set $\gamma^\nu = \lambda^{\bar s}$ and make the GN update~\eqref{GN}.  
The fact that there is a finite value of $s$ for which \eqref{armijo ineq} is satisfied 
follows from inequality $f'(x^\nu; d^\nu) \leq \tilde\Del f(x^\nu) < 0$.
The inequality \eqref{armijo ineq} is called the Armijo inequality. 
A general convergence theory for this algorithm as well as a wide range of others is
found in \cite{Burke85}. 
For the NLLS case, the situation is simple, since $\rho$ is a quadratic, and standard 
convergence theory is given for example in~\cite{Den}. However, the more general theory
is essential in the later sections.

\subsection{Details for Kalman Smoothing}

To implement the Gauss-Newton method described above, one must compute the 
solution $d^\nu$ to the Gauss-Newton subproblem~\eqref{GNsub} for~\eqref{fullNLLS}.
%
%
That is, one must compute
\begin{equation}\label{GNsubNLLS}
d^\nu = \arg\min_{d} \tilde f(d) = \frac{1}{2}\|G^\nu d - \underbrace{ w - g(x^\nu) }_{w^\nu}\|_{Q^{-1}}^2 + \frac{1}{2}\|H^\nu d - \underbrace{ z - h(x^\nu) }_{z^\nu}\|_{R^{-1}}^2\;,
\end{equation}
where 
\begin{equation}\label{NLLSdefs}
G^\nu   = \begin{bmatrix}
    \R{I}  & 0      &          &
    \\
    -g_2^{(1)}(x_1^\nu)   & \R{I}  & \ddots   &
    \\
        & \ddots &  \ddots  & 0
    \\
        &        &   -g_N^{(1)}(x_{N-1}^\nu) & \R{I}
\end{bmatrix}\;, 
\quad
H^\nu = \R{diag}\{h_1^{(1)}(x_1), \dots, h_N^{(1)}(x_N)\}\;.
\end{equation}
However, the problem~\eqref{GNsubNLLS} has exactly the same structure as~\eqref{fullLS}; 
a fact that we have emphasized by defining 
\begin{equation}\label{rhsNU} 
w^\nu := w - g(x^\nu)\;, \quad z^\nu  = z - h(x^\nu)\;.
\end{equation}
Therefore, we can solve it efficiently by using Algorithm~\ref{blockAlg}.

The linearization step in~\eqref{GNsubNLLS} should remind the reader of the EKF. 
Note, however, that 
the Gauss-Newton method is iterative, and we iterate until convergence to a local
minimum of~\eqref{fullNLLS}.  We also linearize along the entire 
state space sequence $x^\nu$ at once in~\eqref{GNsubNLLS}, 
rather than re-linearizing as we make our way through the $x^\nu_k$'s.

\subsection{Numerical Example: Van Der Pol Oscillator}
\label{VanDerPol1}

\begin{figure}\label{fig:vanDerPol}
\begin{center}
{\includegraphics[width=3.5in]{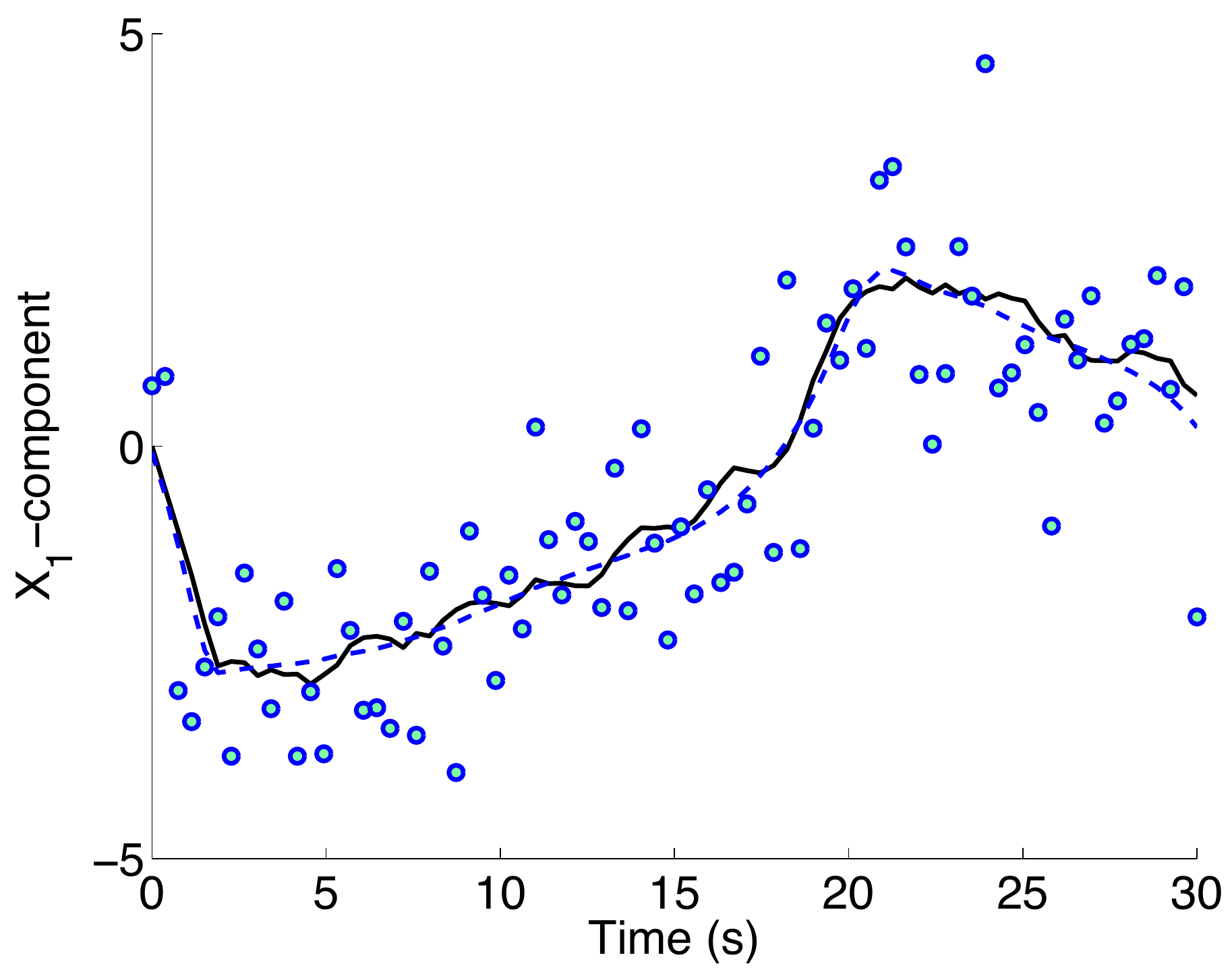}}
{\includegraphics[width=3.5in]{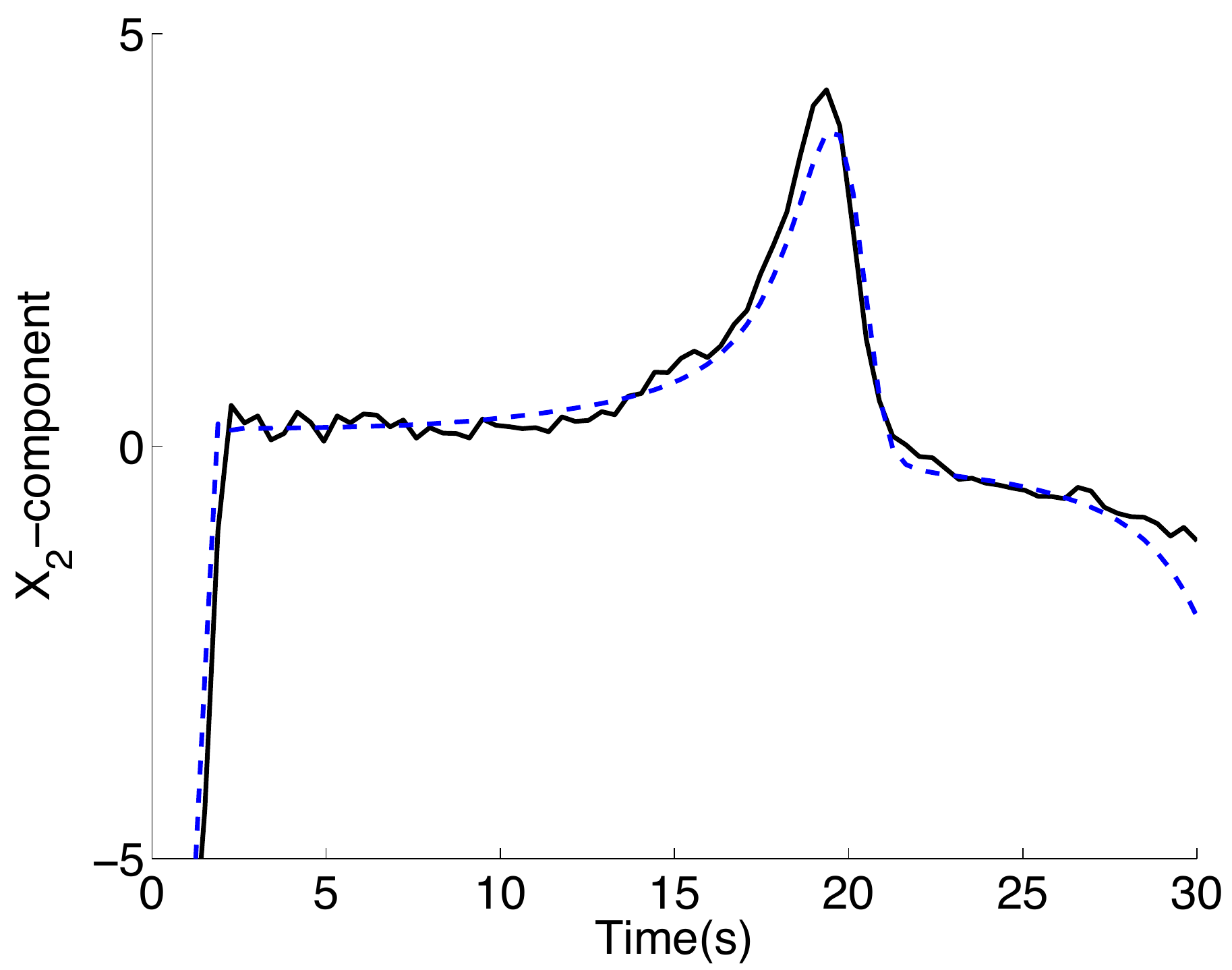}}
\end{center}
\caption{
Tracking the Van Der Pol Osciallator using a nonlinear process model~\eqref{EulerVanODE} and direct (noisy) measurements~\eqref{measEq}
of {\bf  $X_1$-component only}. 
Black solid line is true signal, blue dashed line is nonlinear Kalman smoother estimate. Measurements are displayed as circles. 
}
\end{figure}

The Van der Pol oscillator is a popular
nonlinear process for comparing Kalman filters, see
\cite{GMCDBR06} and \cite[Section 4.1]{KFI08}.
The oscillator is governed by a nonlinear ODE model 
\begin{equation}\label{vanODE}
\dot{X}_1 (t) = X_2 (t)
\hspace{.5cm} \mbox{and} \hspace{.5cm}
\dot{X}_2 (t) = \mu [ 1 - X_1 (t)^2 ] X_2 (t) - X_1 (t)
\; .
\end{equation}
In contrast to the linear model~\eqref{linearSmooth}, which 
was a {\it generic} process for a smooth signal, we now take 
the Euler discretization of~\eqref{vanODE} to be the specific
process model for this situation. 

Given \( X( t_{k-1} ) = x_{k-1} \) the Euler approximation for
\( X ( t_{k-1} + \Delta t ) \) is
\begin{equation}\label{EulerVanODE}
g_k ( x_{k-1} ) = \left( \begin{array}{cc}
	x_{1,k-1} + x_{2,k-1} \Delta t
	\\
	x_{2,k-1} + \{ \mu [ 1 - x_{1,k}^2 ] x_{2,k} - x_{1,k} \} \Delta t
\end{array} \right) \; .
\end{equation}
For the simulation, the `ground truth' is obtained from
a stochastic Euler approximation
of the Van der Pol oscillator.
To be specific,
with \( \mu = 2 \), \( N = 80 \) and \( \Delta t = 30 / N \),
the ground truth state vector \( x_k \) at time \( t_k = k \Delta t \)
is given by \( x_0 = ( 0 , -0.5 )^\R{T} \) and
for \( k = 1, \ldots , N \),
\begin{equation}
\label{VanDerPolTruth}
x_k = g_k ( x_{k-1} ) + w_k \; ,
\end{equation}
where $\{ w_k \}$ is a realization of
independent Gaussian noise with variance $0.01$
and $g_k$ is given in~\eqref{EulerVanODE}.
Our process model for state transitions is also~\eqref{VanDerPolTruth},  
with \( Q_k = 0.01 \; I \) for \( k > 1 \), and so is identical to the model used
to simulate the ground truth \( \{ x_k \} \).
Thus, we have precise knowledge of the process that generated the
ground truth \( \{ x_k \} \).
The initial state \( x_0 \) is imprecisely specified by setting
\( g_1 ( x_0 ) = ( 0.1 , -0.4 )^\R{T} \neq x_0 \)
with corresponding variance \( Q_1 = 0.1 \; I \).
For \( k = 1 , \ldots , N \) noisy measurements \( z_k \) direct measurements 
of {\it the first component only} were used
\begin{equation}\label{measEq}
z_k = x_{1,k} + v_k \; ,
\end{equation}
with $v_k \sim N(0, 1)$. 

The resulting fit is shown in Figure~\ref{fig:vanDerPol}. 
Despite the noisy measurements 
of only $X_1$, we are able to get a good fit for both components. 
The figure was generated using the \verb{ckbs{ package~\cite{ckbs}, see the file \verb{vanderpol_experiment_simple.m{.
The program and example are available for download from COIN-OR.

\section{State space constraints}
\label{Constrained}

In almost every real-world problem, additional prior information is known about the state.  
In many cases, this information can be represented using {\it state space constraints}. 
For example, in tracking physical bodies, we often know (roughly or approximately) 
the topography of the terrain; this information can be encoded as a simple box constraint
on the state. We may also know physical limitations (e.g. maximum acceleration or velocity) 
of objects being tracked, or hard bounds set by biological or financial systems. 
These and many other examples can be formulated using state space constraints.  
The ability to incorporate this information is particularly useful when measurements 
are inaccurate or far between. 

In this section, we first show how to add affine inequality constraints to 
the affine smoother formulation in Section~\ref{Classic}. 
This requires a novel methodology: {\it interior point (IP) methods}, 
an important topic in optimization~\cite{Wright:1997, KMNY91, NN94}. 
IP methods work directly with optimality conditions, so we derive these 
conditions for the smoothing problem. 
Rather than review theoretical results about IP methods, 
we give a general overview and show how they specialize to the linear constrained smoother. 
The constrained Kalman smoother was originally proposed in~\cite{Bell2008}, but we improve 
on that work here, and present a simplified algorithm, which is also faster and more numerically 
stable. We illustrate the algorithm using a numerical example, building on the example in Section~\ref{Classic}.

Once the {\it linear} smoother with {\it linear} inequality constraints 
is understood, we review the constrained {\it nonlinear} smoother
(which can have nonlinear process, measurement, and constraint functions). 
Using~\cite{Bell2008} and references therein, we show that the constrained 
nonlinear smoother is iteratively solved using linear constrained smoothing subproblems,
analogously to how the nonlinear smoother in Section~\ref{Nonlinear}
is iteratively solved using linear smoothing subproblems from Section~\ref{Classic}.  Because 
of this hierarchy, the improvements to the affine algorithm immediately carry over 
to the nonlinear case.  We end with a nonlinear constrained numerical example.  

\subsection{Linear Constrained Formulation}

We start with the linear smoothing problem~\eqref{fullLS}, and 
impose linear inequality constraints on the state space $x$: 
\begin{equation}\label{linearIneq}
B_k x_k \leq b_k\;.
\end{equation}
By choosing the matrix $B_k$ and $b_k$ appropriately, one can ensure $x_k$ lies 
in any polyhedral set, since such a set is defined by a finite intersection of hyperplanes. 
Box constraints, one of the simplest and useful tools for modeling ($l_k \leq x_k \leq u_k$)
can be imposed via 
\[
\begin{bmatrix}I \\ - I \end{bmatrix} x_k \leq \begin{bmatrix} u_k \\ -l_k\end{bmatrix}\;.
\] 
In order to formulate the problem for the entire state space sequence, we define
\begin{equation}\label{linearIneqFull}
B = \R{diag}(\{B_k\})\;, \quad b = \R{vec}(\{b_k\})\;,
\end{equation}
and all of the constraints can be written simultaneously as $Bx \leq b$.
The constrained optimization problem is now given by 
\begin{equation}\label{fullLSc}
\begin{aligned}
\min_{x} f(x) =   &\quad \frac{1}{2}\|Hx - z\|_{R^{-1}}^2 + \frac{1}{2}\|Gx - w\|_{Q^{-1}}^2\\
\text{subject to } &\quad Bx + s = b, \quad s \geq 0\;.
\end{aligned}
\end{equation}
Note that we have rewritten the inequality constraint as an equality constraint 
by introducing a new `slack' variable $s$. 

We derive the Karush-Kuhn-Tucker (KKT) conditions using the Lagrangian formulation. 
The Lagrangian corresponding to~\eqref{linearIneqFull} is given by 
\begin{equation}\label{LagLinearIneq}
\Sc{L}(x, u, s) = \frac{1}{2}\|Hx - z\|_{R^{-1}}^2 + \frac{1}{2}\|Gx - w\|_{Q^{-1}}^2 + u^\top(Bx +s - b)\;.
\end{equation}
The KKT conditions are now obtained by differentiating $\Sc{L}$ with respect to its arguments. 
Recall that the gradient of~\eqref{fullLS} is given by 
\[
(H^\top R^{-1} H + G^\top Q^{-1} G) x  - H^\top R^{-1}z - G^\top Q^{-1}w\;.
\]
As in~\eqref{hessianApprox} set $C = H^\top R^{-1} H + G^\top Q^{-1} G$, 
and for convenience set 
\begin{equation}\label{defGrad}
c =  H^\top R^{-1}z + G^\top Q^{-1}w
\end{equation}
The KKT necessary and sufficient conditions for optimality are given by  
\begin{equation}\label{KKTlin}
\begin{aligned}
\nabla_x \Sc{L} = Cx +c + B^\top u &= 0\\
\nabla_q \Sc{L} = Bx + s - b &= 0\\
 u_i s_i &= 0 \quad \forall i\;; u_i, s_i \geq 0\;.
\end{aligned}
\end{equation}
The last set of nonlinear equations is known as {\it complementarity} conditions. 
In primal-dual interior point methods,
the key idea for solving~\eqref{fullLSc} is to successively solve relaxations of the system
\eqref{KKTlin} that converge to a triplet $(\bar x, \bar u, \bar s)$ which satisfy~\eqref{KKTlin}.

\subsection{Interior Point Approach}
\label{linearIP}
IP methods work directly to find solutions of~\eqref{KKTlin}.  
They do so by iteratively relaxing the complementarity 
conditions $u_i s_i = 0$ to $u_i s_i = \mu$ as they drive 
the relaxation parameter $\mu$ to $0$. The 
{\it relaxed} KKT system is defined by 
\begin{equation}
\label{KKTrel}
F_\mu(s, u, x) 
= 
\begin{bmatrix} 
s  + Bx -b\\
SU\B{1} - \mu\B{1}\\
Cx + B^Tu -c
\end{bmatrix}\;. 
\end{equation}
where $S$ and $U$ are diagonal matrices with $s$ and $u$ on the diagonal, and 
so the second equation in $F_\mu$ implements the relaxation $u_i s_i = \mu$ of~\eqref{KKTlin}. 
Note that the relaxation requires that $\mu_i, s_i > 0$ for all $i$.
Since the solution to~\eqref{fullLSc} is found by driving the KKT system to $0$, 
at every iteration IP methods attempt to drive $F_\mu$ to $0$ 
by Newton's method for root finding. 

Newton's root finding method solves the linear system 
\begin{equation}\label{NewtonKKT}
F_\mu^{(1)}(s, u, x) \begin{bmatrix} \Delta s \\ \Delta u \\ \Delta x\end{bmatrix} = -F_\mu(s, u, x)\;.
\end{equation}
It is important to see the full details of solving~\eqref{NewtonKKT} in order
to see why it is so effective for constrained Kalman smoothing. 
The full system is given by 
\begin{equation}
\label{SimpleNewton}
\begin{bmatrix}
I & 0 & B \\
U & S & 0\\
0 & B^T & C \\
\end{bmatrix}
\begin{bmatrix}
\Delta s \\ \Delta u \\ \Delta x
\end{bmatrix}
= 
-
\begin{bmatrix} 
s  + Bx - b \\
SU\B{1} - \mu\B{1}\\
Cx + B^Tu -c
\end{bmatrix}\;.
\end{equation}
Applying the row operations
\[
\begin{array}{ccc}
\R{row}_2 &\gets& \R{row}_2 - U\R{row}_1 \\
\R{row}_3 &\gets& \R{row}_3 - B^TS^{-1}\R{row}_2
\end{array}\;,
\]
we obtain the equivalent system 
\[
\begin{bmatrix}
I & 0 & B \\
0 & S & -UB\\
0 & 0 & C + B^TS^{-1}UB \\
\end{bmatrix}
\begin{bmatrix}
\Delta s \\ \Delta u \\ \Delta x
\end{bmatrix}
= 
-
\begin{bmatrix} 
s  + Bx - b \\
-U(Bx-b) - \mu\B{1}\\
Cx + B^Tu -c + B^TS^{-1}
\left(U(Bx-b) + \mu \B{1}
\right)
\end{bmatrix}\;.
\]
In order to find the update for $\Delta x$, we have to solve the system 
\begin{equation}\label{xUpdate}
\left(C + B^TS^{-1}UB\right)\Delta x = Cx + B^Tu -c + B^TS^{-1}
\left(U(Bx-b) + \mu \B{1}
\right)
\end{equation}
Note the structure of the matrix in the LHS of~\eqref{xUpdate}. 
The matrix $C$ is the same as in~\eqref{fullLS}, so it is positive definite
symmetric block tridiagonal. The matrices $S^{-1}$ and $U$ are diagonal, 
and we always ensure they have only positive elements. The matrices
$B$ and $B^\top$ are both block diagonal. Therefore, 
$C + B^TS^{-1}UB$ has the same structure as $C$, and 
we can solve~\eqref{xUpdate} using Algorithm~\ref{blockAlg}.

Once we have $\Delta x$, the remaining two updates are obtained by back-solving: 
\begin{equation}\label{uUpdate}
\Delta u = US^{-1}(B(x + \Delta x)-b) + \frac{\mu}{s}
\end{equation}
and
\begin{equation}\label{sUpdate}
\Delta s = -s +b -B(x + \Delta x)\;.
\end{equation}

This approach improves the algorithm presented in~\cite{Bell2008} 
solely by changing the order of variables and equations in~\eqref{KKTrel}.
This approach simplifies the derivation while also improving speed and numerical stability.

It remains to explain how $\mu$ is taken to $0$. There are several strategies,
see~\cite{Wright:1997, KMNY91, NN94}. For the Kalman smoothing application, 
we use one of the simplest: for two out of every three iterations $\mu$ is aggressively 
taken to $0$ by the update $\mu = \mu/10$; while in the remaining iterations, $\mu$ is unchanged. 
In practice, one seldom needs more than 10 interior point iterations; therefore the 
constrained linear smoother performs at a constant multiple of work of the linear smoother.

\begin{figure}\label{fig:affineConstrained}
\begin{center}
{\includegraphics[width=3.5in]{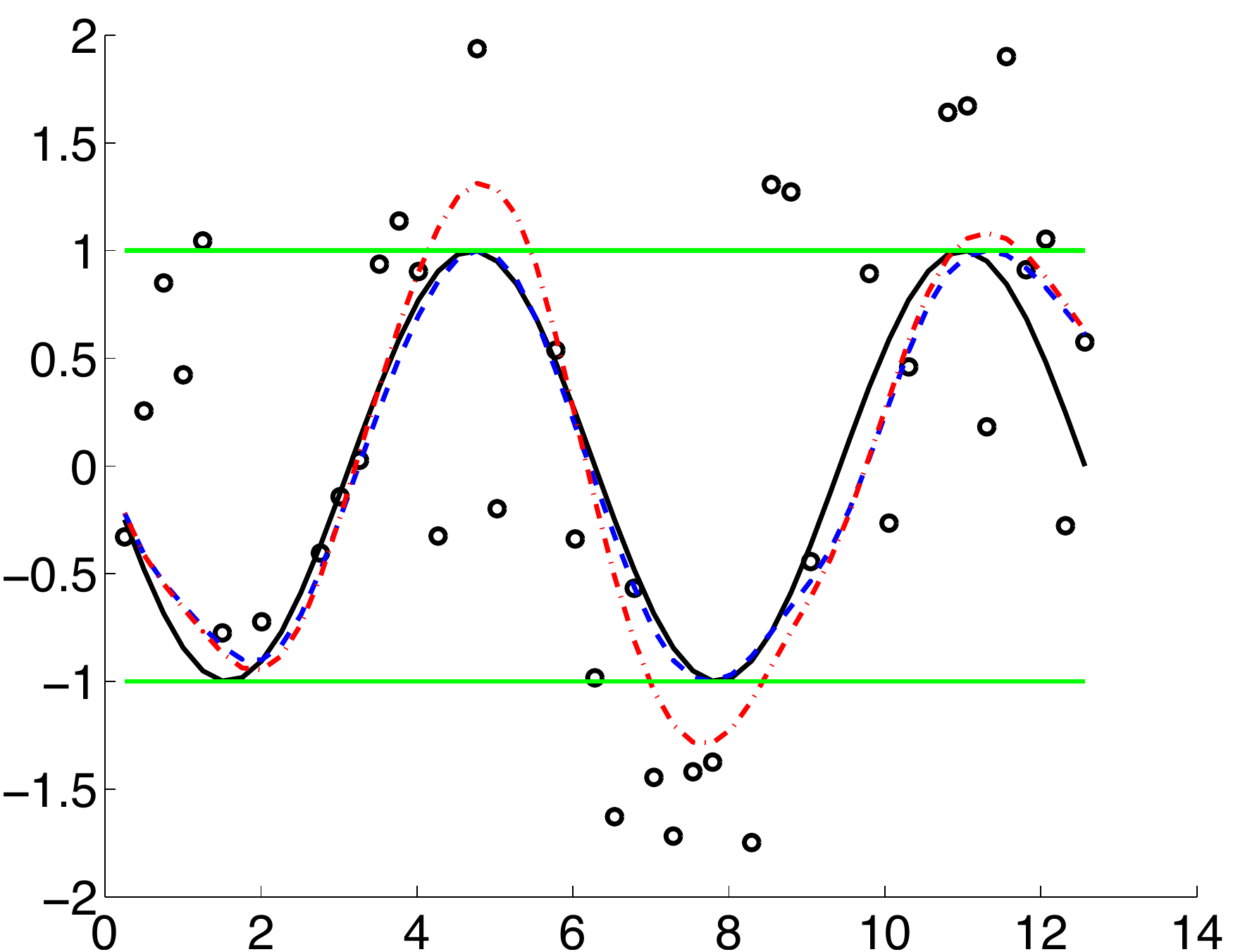}}
{\includegraphics[width=3.5in]{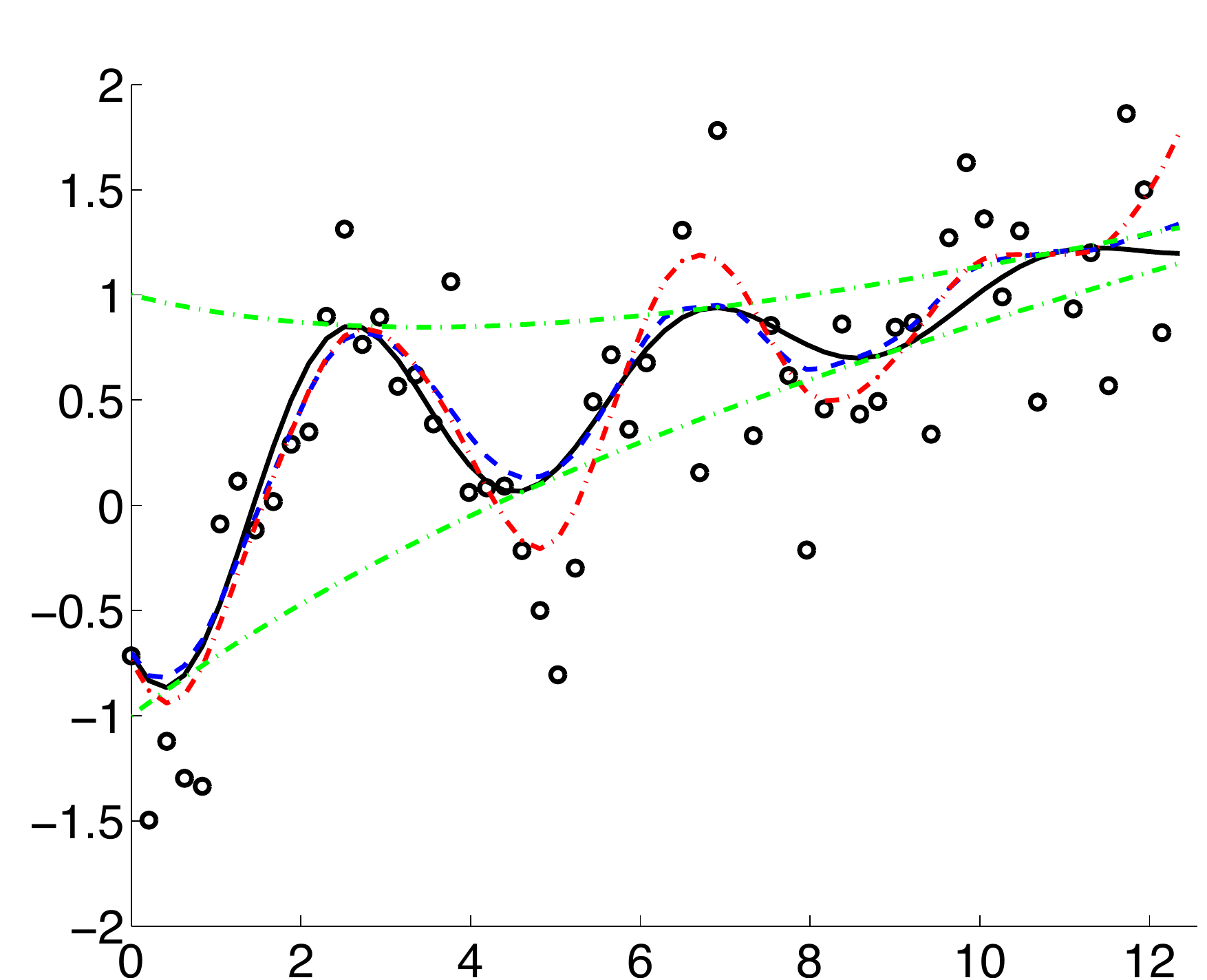}}
\end{center}
\caption{
Two examples of linear constraints. 
Black solid line is true signal, magenta dash-dot line is unconstrained Kalman smoother, and blue dashed line 
is the constrained Kalman smoother. Measurements are displayed as circles, and bounds are shown as 
green horizontal lines. In the left panel, note that performance of the bounded smoother 
is significantly better around time 4-10 --- the unconstrained is fooled by the  measurements at times 4 and 8. 
In the right panel, as the oscillations die down due to damping, the measurement variance remains unchanged, 
so it becomes much more difficult to track the signal without the bound constraints. 
}
\end{figure}

\subsection{Two Linear Numerical Examples }

In this section, we present to simple examples, both with linear constraints. 

\paragraph{Constant Box Constraints}
In the first example, we impose box constraints in the example of Section~\ref{sec:LinearExample}.
Specifically, we take advantage of the fact the state is bounded: 
\begin{equation}\label{boxSin}
\begin{bmatrix} -1\end{bmatrix} \leq \begin{bmatrix}x \end{bmatrix} \begin{bmatrix} 1 \end{bmatrix}
\end{equation}

We can encode this information in form~\eqref{linearIneq} with 
\begin{equation}\label{boxSinFull}
B_k = \begin{bmatrix} 1 & 0 \\ 0 & -1  \end{bmatrix}\;, \quad b_k = \begin{bmatrix} 1\\1 \end{bmatrix}\;.
\end{equation}

We contrast the performance of the constrained linear smoother with that of the linear smoother without constraints. 
To show the advantages of modeling with constraints, we increase the measurement noise 
in both situations to $\sigma^2 = 1$. The results are show in Figure~\ref{fig:affineConstrained}. 
The constrained smoother avoids some of the problems encountered by the unconstrained smoother. 
Of particular interest are the middle and end parts of the track, where the unconstrained smoother goes far afield because
of bad measurement. The constrained smoother is able track portions of the track extremely well, having 
avoided the bad measurements with the aid of the bound constraints. The figure was generated 
using the \verb{ckbs{ package~\cite{ckbs}, specifically using the example file \verb{affine_ok_boxC.m{.

\paragraph{Variable Box Constraints}

In the second example, we impose time-varying constraints on the state. Specifically, we track an exponentially bounded
signal with a linear trend: 
\[
\exp(-\alpha t)\sin(\beta t) + .1t
\]
using the `smooth signal' process model and direct measurements, as in Section~\ref{sec:LinearExample}. The challenge here is that 
as the oscillations start to die down because of the exponential damping, the variance of the measurements
remains the same. We can improve the performance by giving the smoother the exponential damping terms as constraints. 

We included the second example to emphasize that `linearity' of constraints means `with respect to the state'; in fact, 
the constraints in the second example are simply box constraints which are time dependent. The second example 
is no more complicated than the first one for the constrained smoother.


\subsection{Nonlinear Constrained Smoother}
\label{sec:ConvexCompositeCon}

We now consider the nonlinear constrained smoother, where we allow process functions $g_k$, measurement functions $h_k$
to be nonlinear, and also allow nonlinear smooth constraints $\xi_k(x_k) \leq b_k$. To be consistent with the notation we use throughout 
the paper, we define a new function
\begin{equation}\label{fullCon}
\xi(x) = \begin{bmatrix} \xi_1(x_1) \\ \xi_2(x_2)\\ \vdots \\ \xi_N(x_N) \end{bmatrix}\;,
\end{equation}
so that all the constraints can be written simultaneously as $\xi(x) \leq b$.

The problem we would like to solve now is a constrained reformulation of~\eqref{fullNLLS}
\begin{equation}
\label{conNLLS}
\begin{aligned}
\min_x \quad &f(x) = \frac{1}{2}\|g(x) - w\|_{Q^{-1}}^2 + \frac{1}{2}\|h(x) - z\|_{R^{-1}}^2 \\
\text{subject to } \quad & \xi(x) -b\leq 0\;.
\end{aligned}
\end{equation}

At this point, we come back to the convex-composite representation described in Section~\ref{sec:ConvexComposite}. 
The constraint $\xi(x) - b \leq 0$ may be represented using an additional term in the objective function: 
\begin{equation}\label{conIndicator}
\indicator{\xi(x) - b}{\mathbb{R}_{-}}\;,
\end{equation}
where $\indicator{x}{C}$ is the convex indicator function: 
\begin{equation}\label{indicator}
\indicator{x}{C} = \begin{cases}   0 & x \in C \\ \infty &x \not \in C\end{cases}\;.
\end{equation}

Therefore, the objective~\eqref{conNLLS} can be represented as follows: 

\begin{equation}\label{convexCompositeCon}
\begin{aligned}
f(x) &= \rho(F(x))\\
\rho\left(\begin{matrix} y_1\\ y_2\\ y_3\end{matrix}
\right) 
&= \frac{1}{2}\|y_1\|_{Q^{-1}}^2 + \frac{1}{2}\|y_2\|_{R^{-1}}^2 + \indicator{y_3}{\mathbb{R}_-}\\
\quad F(x) &= \begin{bmatrix} g(x) - w \\ h(x) - z \\ \xi(x) - b\end{bmatrix}\;.
\end{aligned}
\end{equation}

The approach to nonlinear smoothing in~\cite{Bell2008} is essentially the Gauss-Newton method 
described in Section~\ref{sec:GaussNewton}, applied to~\eqref{convexCompositeCon}. 
In other words, at each iteration $\nu$, the function $F$ is linearized, and the direction finding 
subproblem is obtained by solving 

\begin{equation}\label{conGNsub}
\begin{aligned}
\min_d  & \quad \frac{1}{2}\|G^\nu d - \underbrace{ w - g(x^\nu) }_{w^\nu}\|_{Q^{-1}}^2 + \frac{1}{2}\|H^\nu d - \underbrace{ z - h(x^\nu) }_{z^\nu}\|_{R^{-1}}^2\;, \\
\text{subject to } &  \quad B^\nu d \leq \underbrace{b - \xi(x^\nu)}_{b^\nu}\;,
\end{aligned}
\end{equation}
where $G^\nu$ and $H^\nu$ are exactly as in~\eqref{GNsubNLLS}, $B^\nu = \nabla_x \xi(x^\nu)$ is a block diagonal matrix because of the structure of $\xi$~\eqref{fullCon},
and we have written the indicator function in~\eqref{convexCompositeCon} as an explicit constraint to emphasize the structure of the subproblem.  

Note that~\eqref{conGNsub} has exactly the same structure as the linear constrained smoothing problem~\eqref{fullLSc}, and therefore can be solved using 
the interior point approach in the previous section.  Because the convex-composite objective~\eqref{convexCompositeCon} is not finite valued
(due to the indicator function of the feasible set), to prove convergence of the nonlinear smoother,~\cite{Bell2008} uses results from~\cite{BurkeHan1989}.
We refer the interested reader to~\cite[Lemma 8, Theorem 9]{Bell2008} for theoretical convergence results, and to~\cite[Algorithm 6]{Bell2008} for 
the full algorithm, including line search details. 

Because of the hierarchical dependence of the nonlinear constrained smoother on the linear constrained smoother, the simplified improved approach
we presented in Section~\ref{linearIP} pays off even more in the nonlinear case, where it is used repeatedly as a subroutine. 

\subsection{Nonlinear Constrained Example}
\label{NonlinearConstrained}

\begin{figure}\label{fig:sineConstrained}
\begin{center}
{\includegraphics[width=4in]{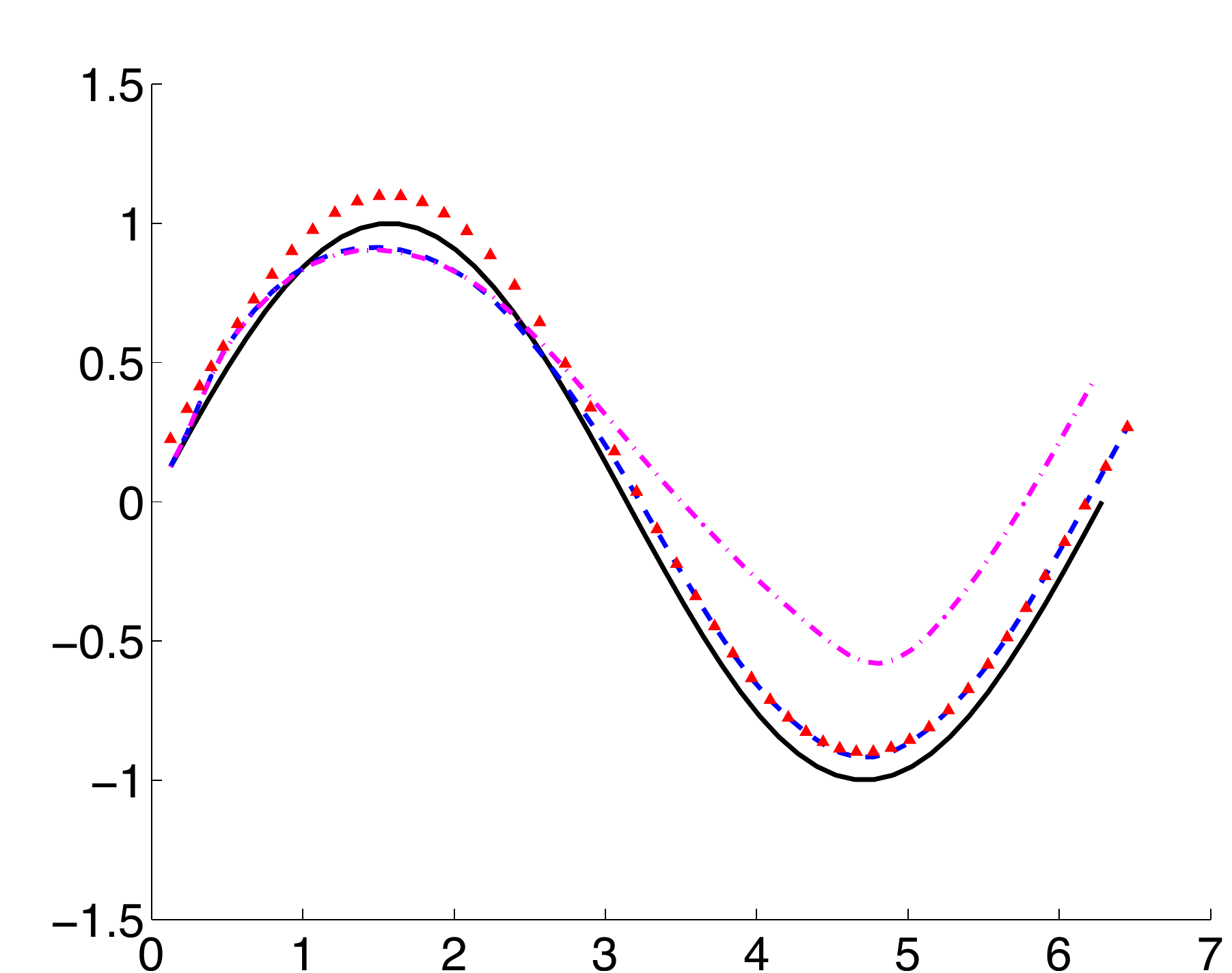}}
\end{center}
\caption{
Smoother results for ship tracking example with linear process model, nonlinear measurement model, and nonlinear constraints (with respect to the state). 
Black solid line is true state, red triangles denote the constraint, magenta dash-dot line is the unconstrained estimate, 
and blue dashed line gives the constrained nonlinear smoothed estimate. 
}
\end{figure}

The example in this section is reproduced from~\cite{Bell2008}. 
Consider the problem of tracking a ship traveling
close to shore where we are given distance measurements from two
fixed stations to the ship as well as the location of the shoreline.
Distance to fixed stations is a nonlinear function, so the measurement model
here is nonlinear. 
 
In addition, the corresponding constraint functions \( \{ f_k \} \)
are not affine because the shoreline is not a straight line.
For the purpose of simulating the measurements \( \{ z_k \} \),
the ship velocity \( [ X_1 (t) , X_3 (t) ] \)
and the ship position \( [ X_2 (t) , X_4 (t) ] \)
are given by 
\[
X(t) = [ ~ 1 ~ , ~ t ~ , ~ -\cos (t) ~ , ~ 1.3 - \sin (t) ~ ]^\top
\]
Both components of the ship's position are modeled using the smooth signal 
model in Section~\ref{sec:LinearExample}. Therefore we introduce two velocity
components, and the process model is given by 
\[G_k
= 
\begin{bmatrix} 
1 & 0 & 0 & 0 \\
\Delta t & 1 & 0 & 0 \\
0 & 0 & 1 & 0\\
0 & 0 & \Delta t & 0
\end{bmatrix}\;,  \quad
Q_k =\begin{bmatrix}
    \Delta t       & \Delta t^2 / 2  & 0              & 0              \\
    \Delta t^2 / 2 & \Delta t^3 / 3  & 0              & 0              \\
    0              & 0               & \Delta t       & \Delta t^2 / 2 \\
    0              & 0               & \Delta t^2 / 2 & \Delta t^3 / 3
\end{bmatrix}\;.
\]
The initial state estimate is given by
\(
g_1 ( x_0 ) = X( t_1 )
\)
and
\(
Q_1 =  100 I_4
\)
where \( I_4 \) is the four by four identity matrix.
The measurement variance
is constant for this example and is denoted by \( \sigma^2 \).
The distance measurements are made from two stationary locations on shore.
One is located at $(0, 0)$ and the other is located at  $(2 \pi , 0 )$.
The measurement model is given by 
\[
h_k ( x_k ) = \left( \begin{array}{c}
    \sqrt{ x_{2,k}^2 + x_{4,k}^2 }
    \\
    \sqrt{ ( x_{2,k} - 2 \pi )^2 + x_{4,k}^2 }
\end{array} \right)
\; , \;
R_k = \left( \begin{array}{cc}
    \sigma^2 & 0 \\
    0        & \sigma^2
\end{array} \right)\;.
\]

We know that the ship does not cross land, so
\(
    X_4 (t) \geq 1.25 - \sin [ X_2 (t) ]
\).
This information is encoded by the constraints
\[
\xi_k ( x_k ) = 1.25 - \sin ( x_{2,k} ) - x_{4,k} \leq 0\;.
\]
The initial point for the smoother is $[0, 0, 0, 1]^\top$,which is not feasible. The results are plotted in Figure~\ref{fig:sineConstrained}. 
The constrained smoother performs significantly better than the unconstrained smoother in this example. The experiment was done
using the \verb{ckbs{ program, specifically see \verb{sine_wave_example.m{.

\section{Robust Kalman smoothing}
\label{Robust}  

In many applications, the probalistic model for the dynamics and/or the
observations \eqref{IntroGaussModel} is not well described by a Gaussian
distribution. This occurs in the model for the
observations when they are contaminated by outliers, or more generally, 
when the measurement noise
$v_k$ is heavy tailed \cite{Schick1994}, and it occurs in the model for the dynamics
when tracking systems with rapidly changing dynamics, or jumps in the state values \cite{Boyd2009}.
A {\em robust} Kalman filter or smoother is one that can obtain an acceptable estimate of the
state when Gaussian assumptions are violated, and which continues to perform well when they
are not violated. 

We show how to accommodate non-Gaussian densities by starting with a simple case of
non-Gaussian heavy tailed measurement noise $v_k$ \cite{Aravkin2011tac}.
However, this general approach can be extended to $w_k$ as well.
Heavy tailed measurement noise 
occurs in applications
related to glint noise \cite{HewerMartin1987}, 
turbulence, asset returns, and sensor failure or machine malfunction. 
It can also occur in the presence of secondary noise sources or other kinds of data anomalies.
Although it is possible to estimate a minimum variance estimate of the state using
stochastic simulation methods such as Markov chain Monte-Carlo (MCMC) or particle filters
\cite{GMCDBR06,Liu98}, these methods are very computationally intensive, and convergence often relies
on heuristic techniques and is highly variable. 
The approach taken here is very different. It is based on the
optimization perspective presented in the previous sections. We develop a method for
computing the MAP estimate of the state sequence under the assumption that the
observation noise comes from the $\ell_1$-Laplace density often used in
robust estimation, e.g., see \cite[equation 2.3]{Gao2008}. 
As we will see, the resulting optimization problem will again be one of convex composite 
type allowing us to apply a Gauss-Newton strategy for computing the MAP estimate. Again,
the key to a successful computational strategy is the preservation of the underlying 
tri-diagonal structure.

\subsection{An $\ell_1$-Laplace Smoother}

\begin{figure*}
\center
{\includegraphics[scale=0.4]{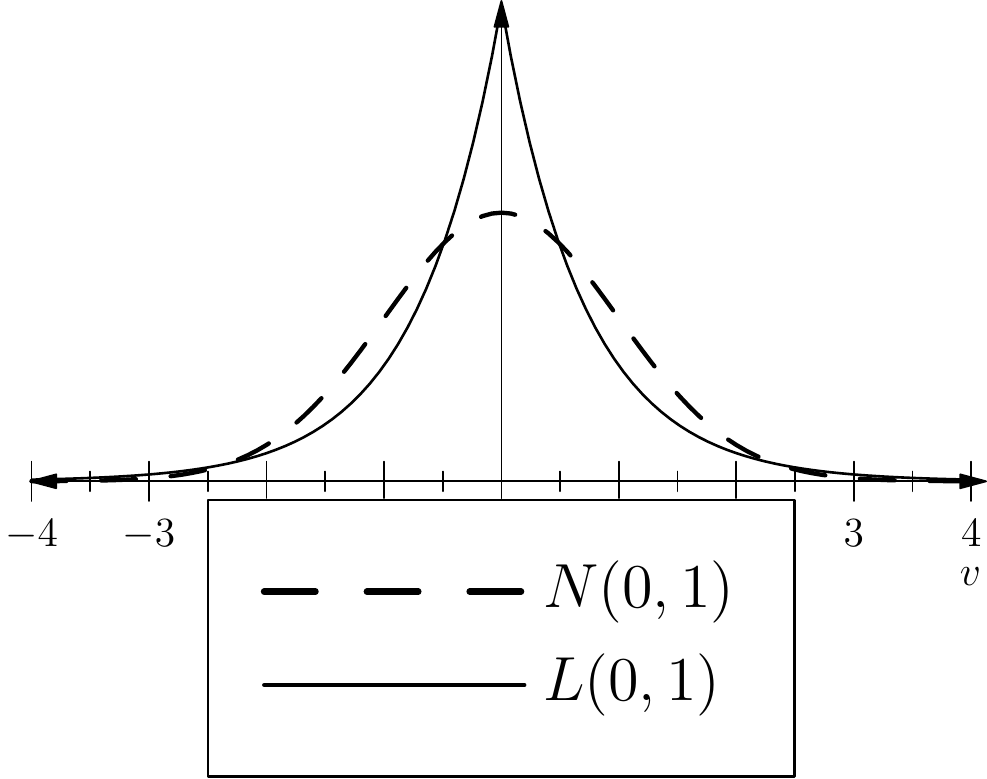}}{\includegraphics[scale=0.4]{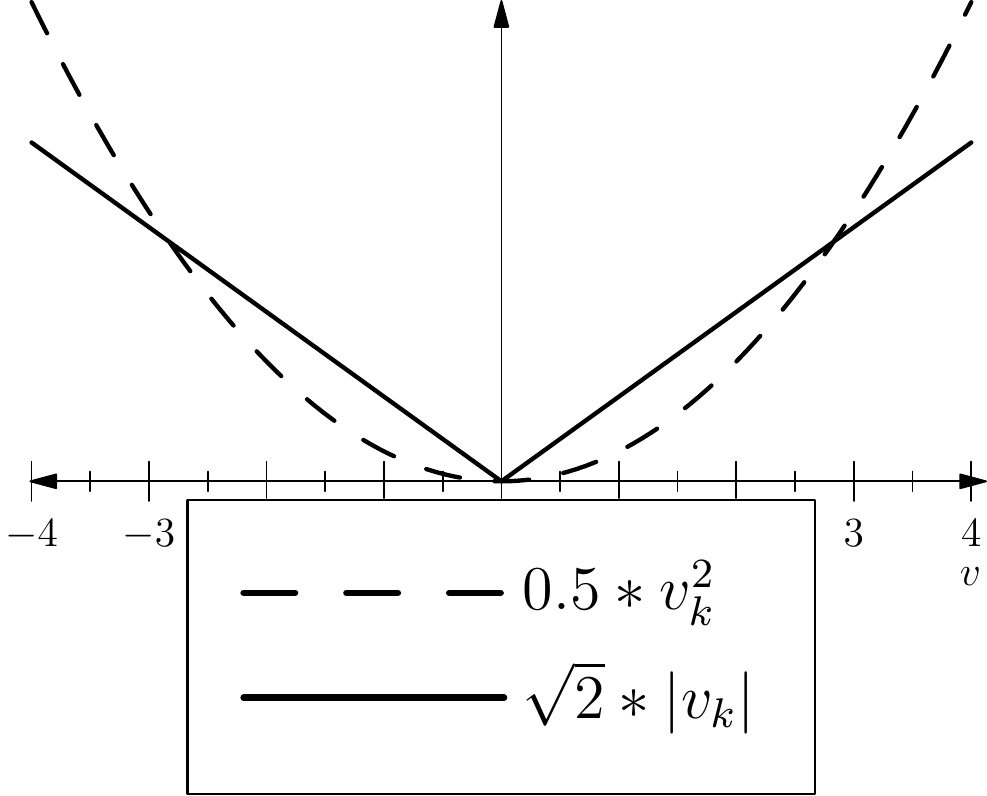}}{\includegraphics[scale=0.4]{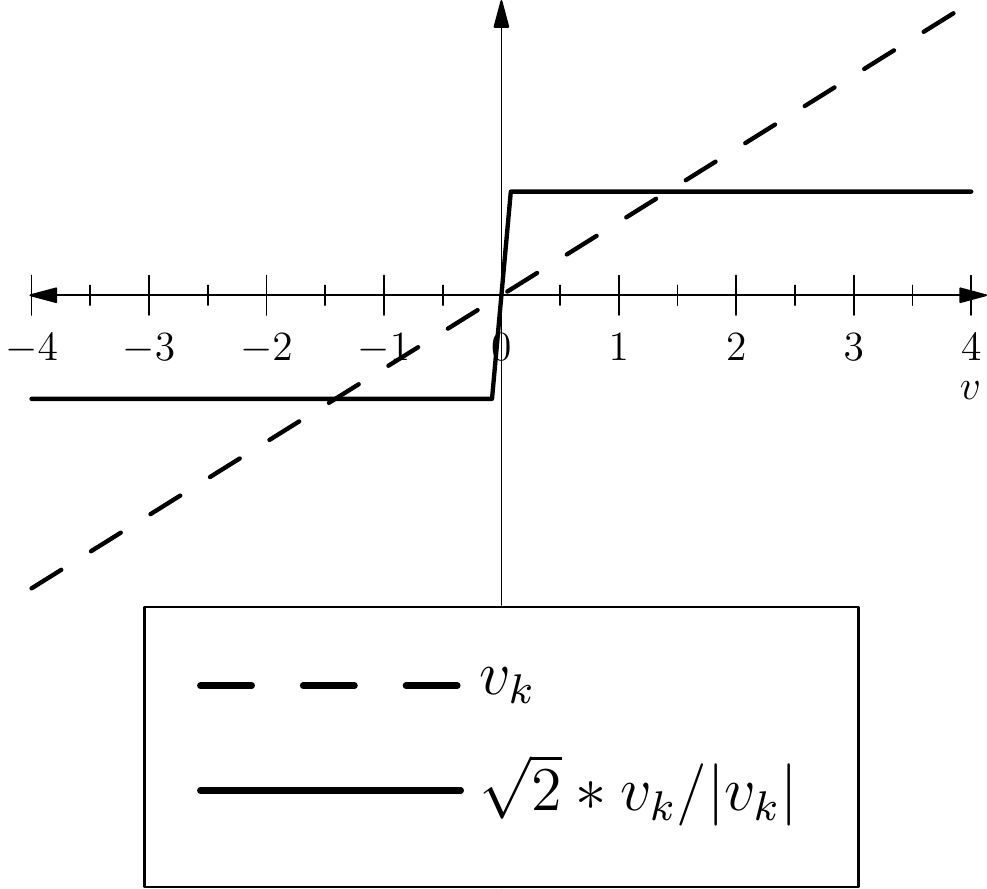}}
\caption{\label{GaussianAndLaplace}
Gaussian and Laplace
Densities, Negative Log Densities, and Influence Functions
(for scalar $v_k$)
}
\end{figure*}

For $u \in \mB{R}^m$ we use the notation
$ \| u \|_1 $ for the $ \ell_1 $ norm of $u$; i.e.,
$\| u \|_1 = | u_1 | + \ldots + | u_m |$.
The multivariate $\ell_1$-Laplace distribution
with mean $\mu$ and covariance $R$
has the following density:
\begin{eqnarray}
\label{LaplaceDensity}
\B{p} ( v_k )
& = &
\det \left( 2 R \right)^{-1/2}
\exp \left[ - \sqrt{2} \left\| R^{-1/2} ( v_k - \mu ) \right\|_1 \; \right] \; ,
\end{eqnarray}
where
$R^{1/2}$ denotes a Cholesky factor of the positive definite matrix $R$;
i.e., $R^{1/2} ( R^{1/2} )^\R{T} = R$.
One can verify that this is a probability distribution with covariance
$R$ using the change of variables $u = R^{-1/2} ( v_k - \mu )$.
A comparison of the Gaussian and Laplace distributions
is displayed in Figure~\ref{GaussianAndLaplace}.
This comparison includes the densities,
negative log densities,
and influence functions,
for both distributions.

\subsubsection{Maximum {\it a posteriori}  formulation}
\label{sec:l1MAP}

Assume that the model for the dynamics and the observations is given by \eqref{IntroGaussModel},
where $w_k$ is assumed to be Gaussian and $v_k$ is modeled by the 
$\ell_1$-Laplace density \eqref{LaplaceDensity}. 
Under these assumptions, the MAP objective function is given by
\begin{equation}\label{l1MAP}
\begin{aligned}
P\left(\{x_k\} \big| \{z_k\}\right) &\propto P\left(\{z_k\}\big|\{x_k\}\right)P\left(\{x_k\}\right)\\
& = \prod_{k=1}^N P(\{v_k\})  P(\{w_k\})\\
& \propto \prod_{k=1}^N \exp\left( - \sqrt{2} \left\| R^{-1/2} ( z_k - h_k( x_k) ) \right\|_1
-\frac{1}{2}(x_k - g_k( x_{k-1}))^\top Q_k^{-1}(x_k - g_k (x_{k-1}))\right)\;.
\end{aligned}
\end{equation}
Dropping terms that do not depend on $ \{ x_k \} $,
minimizing this MAP objective with respect to $ \{ x_k \} $
is equivalent to minimizing
\begin{gather*}
f(\{x_k\}):=\sqrt{2} \sum_{k=1}^N \left\| R_k^{-1/2} [ z_k - h_k (x_k) ] \right\|_1
+ 
\frac{1}{2} \sum_{k=1}^N
	[ x_k - g_k ( x_{k-1} ) ]^\R{T} Q_k^{-1} [ x_k - g_k ( x_{k-1} ) ] \; ,
\end{gather*}
where, as in \eqref{IntroGaussModel}, $x_0$ is known and $g_0=g_1(x_0)$.
Setting
\begin{equation}\label{re-defs}
\begin{aligned}
R       & =  \R{diag} ( \{ R_k \} )
\\
Q       & =  \R{diag} ( \{ Q_k \} )
\\
x       & = \R{vec} ( \{ x_k \} )
\\
w      &  = \R{vec} (\{g_0, 0, \dots, 0\})
\\
z      & = \R{vec} (\{z_1,  z_2, \dots, z_N\})
\end{aligned} ,\quad \quad
g(x) = \begin{bmatrix}x_1 \\ x_2 - g_2(x_1) \\ \vdots \\ x_N - g_N(x_{N-1}) \end{bmatrix}\;,
\quad
h(x) = \begin{bmatrix} h_1(x_1) \\ h_2(x_2) \\ \vdots \\ h_N(x_N)\end{bmatrix}\;,
\end{equation}
as in \eqref{defs} and \eqref{ghDef}, 
%
the MAP estimation problem is equivalent to
\begin{equation}
\label{l1MAP2}
\begin{array}{c}
\mbox{minimize} \\ \displaystyle{ x \in \B{R}^{N n} }
\end{array}
\;
f(x)=	\frac{1}{2}\left\| 	g ( x )-w\right\|_{Q^{-1}}
		+
	\sqrt{2} \left\| R^{-1/2} (h(x)-z) \right\|_1.
\end{equation}

\subsubsection{The Convex Composite Structure}\label{sec:cvxmp3}

The objective in~\eqref{l1MAP2} can again be written as a 
the composition of a convex function $\rho$ with a smooth function $F$:
\begin{equation}\label{convexComposite3}
f(x) = \rho(F(x))\;, 
\end{equation}
where
\begin{equation}
\label{CCparts3}
\rho\left(\begin{matrix} y_1\\ y_2\end{matrix}
\right) 
= \frac{1}{2}\|y_1\|_{Q^{-1}}^2 + \sqrt{2}\|R^{-1/2}y_2\|_1\;, \quad F(x) = \begin{bmatrix} g(x)-w  \\ h(x) - z\end{bmatrix}\;.
\end{equation}
Consequently, the generalized Gauss-Newton methodology described in Section \ref{sec:GaussNewton} 
again applies. 
That is, given an approximate solution $x^\nu$ to \eqref{l1MAP2}, 
we compute a new approximate solution of the form
\[
x^{\nu+1}=x^\nu+\gamma^\nu d^\nu,
\]
where $d^\nu$ solves the subproblem
\begin{equation}\label{cvxcmp2}
\mathop{\mbox{minimize}}_{d\in\Rn}\rho(F(x^\nu)+F'(x^\nu)d),
\end{equation}
and $\gamma^\nu$ is computed using the backtracking line-search procedure described
in Section \ref{sec:GaussNewton}. Following the pattern described in \eqref{GNsubNLLS},
the subproblem \eqref{cvxcmp2}, where $\rho$ and $F$ are given in \eqref{CCparts3},
has the form
\begin{equation}\label{GNsubRobust}
d^\nu = \arg\min_{d} \tilde f(d) = 
\frac{1}{2}\|G^\nu d - \underbrace{ w - g(x^\nu) }_{w^\nu}\|_{Q^{-1}}^2 + 
\sqrt{2}\| R^{-1/2}(H^\nu d - \underbrace{ z - h(x^\nu) }_{z^\nu})\|_1\;,
\end{equation}
where 
\begin{equation}\label{Robustdefs}
G^\nu   = \begin{bmatrix}
    \R{I}  & 0      &          &
    \\
    -g_2^{(1)}(x_1^\nu)   & \R{I}  & \ddots   &
    \\
        & \ddots &  \ddots  & 0
    \\
        &        &   -g_N^{(1)}(x_{N-1}^\nu) & \R{I}
\end{bmatrix}\;, 
\quad
H^\nu = \R{diag}\{h_1^{(1)}(x_1), \dots, h_N^{(1)}(x_N)\}\;.
\end{equation}

\subsubsection{Solving the Subproblem by Interior Point Methods}

By \eqref{GNsubRobust}, the basic subproblem that must be solved
takes the form
\begin{equation}\label{RobustBasicSubproblem}
\min_{d} 
\frac{1}{2}\|G d - w\|_{Q^{-1}}^2 + 
\sqrt{2}\| R^{-1/2}(H d - z)\|_1\;,
\end{equation}
where, as in \eqref{defs},
\begin{equation}\label{Robustdefs}
\begin{aligned}
R       & =  \R{diag} ( \{ R_k \} )
\\
Q       & =  \R{diag} ( \{ Q_k \} )
\\
H       & = \R{diag} (\{H_k\} )
\end{aligned}\quad \quad
\begin{aligned}
x       & = \R{vec} ( \{ x_k \} )
\\
w      &  = \R{vec} (\{w_1,w_2,\dots,w_N\})
\\
z      & = \R{vec} (\{z_1,  z_2, \dots, z_N\})
\end{aligned} \quad \quad
\begin{aligned}
G  & = \begin{bmatrix}
    \R{I}  & 0      &          &
    \\
    -G_2   & \R{I}  & \ddots   &
    \\
        & \ddots &  \ddots  & 0
    \\
        &        &   -G_N  & \R{I}
\end{bmatrix}\;.
\end{aligned}
\end{equation}
Using standard optimization techniques, one can introduce a pair
of auxiliary non-negative variables $p^+,p^-\in \mB{R}^M$
($M=\sum_{k=1}^Nm(k)$) so that this problem can be rewritten as
\begin{equation}\label{Robustsub2}
\begin{array}{ll}
\mbox{minimize}&\frac{1}{2}d^\top Cd+c^\top d +\B{\sqrt{2}}^\top(p^++p^-)\\
\mbox{w.r.t.}&d\in\mB{R}^{nN},\ p^+,p^-\in\mB{R}^M\\
\mbox{subject to}&Bd+b=p^+-p^-\ ,
\end{array}
\end{equation}
where
\[
C =G^\top Q^{-1}G=
\begin{bmatrix}
	C_1  & A^\top_2  & 0 &
	\\
	A_2 & C_2  & A^\top_3  & 0
	\\
	0 & \ddots & \ddots& \ddots &
	\\
	& 0 & A_N  & C_N 
\end{bmatrix}
\; ,\quad\quad
\begin{aligned}
A_k  &= -Q_k^{-1}G_{k} 
\\
C_k  &=
Q_k^{-1} + G_{k+1}^\top Q^{-1}_{k+1} G_{k+1} 
\\
c&= G^\top w
\\
B&=R^{-1/2}H
\\
b&=-R^{-1/2}z
\end{aligned}
\ .
\]
The problem \eqref{Robustsub2} is a convex quadratic program. If we define 
\begin{equation}
\label{DefineFRobust}
F_\mu (
	p^+ ,
	p^- ,
	s^+ ,
	s^-,
	d)
=
\begin{bmatrix}
p^+ - p^- - b - Bd
\\
\R{diag}(p^-) \R{diag}(s^-) \B{1} - \mu \B{1}
\\
s^+ + s^- - 2 \sqrt{\B{2}}
\\
\R{diag}(p^+) \R{diag}(s^+) \B{1} - \mu \B{1}
\\
C d + c + B^\R{T}( s^- - s^+ ) / 2
\end{bmatrix}
\; ,
\end{equation}
for $\mu\ge 0$, then the KKT conditions for \eqref{Robustsub2} 
can be written as
\[
F_0(p^+,p^-,s^+,s^-,d)=0\ .
\]
The set of solutions to $F_\mu(p^+,p^-,s^+,s^-,d)=0$ for $\mu>0$ is called the central path.
We solve the system for $\mu=0$ by an interior point strategy which, as described earlier, is a Newton based predictor-corrector 
method for following the central path as $\mu\downarrow 0$. 
At each iteration of the interior point method we need to solve a system of the form
\[
F_\mu(p^+,p^-,s^+,s^-,d)+F'_\mu(p^+,p^-,s^+,s^-,d)
\begin{bmatrix}
	\Delta p^+ \\
	\Delta p^- \\
	\Delta s^+ \\
	\Delta s^- \\
	\Delta y
\end{bmatrix}
=
0,
\]
where the vectors $ p^+,p^-,s^+,$ and $s^-$ are componentwise strictly positive.
Using standard methods of Gaussian elimination (as in Section \ref{linearIP}),
we obtain the solution
\begin{eqnarray*}
\Delta y\ \,&= &
	[ C + B^\R{T}  T^{-1} B ]^{-1} ( \bar{e} + B^\R{T} T^{-1} \bar{f} )
\\
\Delta s^- &= &
	T^{-1} B \Delta y - T^{-1} \bar{f}
\\
\Delta s^+ &= &
	- \Delta s^- +  2 \sqrt{\B{2}} - s^+ - s^-
\\
\Delta p^- &= &
	\R{diag}( s^- )^{-1} [  \tau \B{1} - \R{diag}( p^- ) \Delta s^- ] - p^-
\\
\Delta p^+ &= &
	\Delta p^- + B \Delta y + b + B y - p^+ + p^-,
\end{eqnarray*}
where
\begin{eqnarray*}
\bar{d}  &= &  \tau \B{1} / s^+  - \tau \B{1} / s^- - b - B y + p^+
\\
\bar{e}  &= & B^\R{T} ( \sqrt{\B{2}} - s^- ) - C y - c
\\
\bar{f}  &= & \bar{d} - \R{diag}( s^+ )^{-1} \R{diag}( p^+ ) ( 2 \sqrt{\B{2}} - s^- )
\\
T        &= &\R{diag}( s^+ )^{-1} \R{diag}( p^+ ) + \R{diag}( s^- )^{-1} \R{diag}( p^- )\ .
\end{eqnarray*}
Since the matrices $T$ and $B$ are block diagonal, the matrix $B^\top TB$ is also block
diagonal.
Consequently, the key matrix $C + B^\R{T}  T^{-1} B$ has exactly the same form as the
block tri-diagonal matrix in \eqref{BlockTridiagonalEquation} with
\begin{equation*}
\begin{aligned}
c_k
&=Q_k^{-1} + G_{k+1}^\top Q^{-1}_{k+1} G_{k+1} +
H_k^\top T_k^{-1}H_k
& k=1,\dots,N,
\\
a_k&=-Q_k^{-1}G_k& k=2,\dots,N,
\end{aligned}
\end{equation*}
where $T_k=\R{diag}( s_k^+ )^{-1} \R{diag}( p_k^+ ) + \R{diag}( s_k^- )^{-1} \R{diag}( p_k^- )$.
Algorithm \ref{blockAlg} can be applied to solve this system accurately and stably with
$O(n^3N)$ floating point operations which preserves the efficiency of the classical Kalman Filter algorithm.

Further discussion on how to incorporate approximate solutions to the quadratic programming
subproblems can be found in \cite[Section V]{Aravkin2011tac}.

\subsubsection{A Linear Example}

In the linear case, the functions $g_k$ and $h_k$ is \eqref{IntroGaussModel} are affine
so that they equal their linearizations. In this case, the problems \eqref{l1MAP2} and
\eqref{cvxcmp2} are equivalent and only one subproblem of the form
\eqref{RobustBasicSubproblem}, or equivalently \eqref{Robustsub2}, 
needs to be solved. We illustrate the $\ell_1$-Laplace
smoother described in Section \ref{sec:l1MAP} by applying it to the example studied
in Section \ref{sec:LinearExample}, except now the noise term $v_k$ is modeled using the
$\ell_1$-Laplace density. The numerical experiment described below
is take from \cite[Section VI]{Aravkin2011tac}.

\begin{table}
\caption{
\label{SimulationResults}
Median MSE and 95\% confidence intervals
for the different estimation methods
}
\begin{center}
\tabcolsep=0.15cm
\begin{tabular}{|c|c|c|c|c|}\hline
$p$&$\phi$
	& GKF
	& IGS
	& ILS
\\ \hline
$0$&$-$
    &.34 (.24, .47)
    &.04(.02, .1)
    &.04(.01, .1)
\\ \hline
$.1$&$1$
    &.41(.26, .60)
    &.06(.02, .12)
    &.04(.02, .10)
\\ \hline
$.1$&$4$
    &.59(.32, 1.1)
    &.09(.04, .29)
    &.05(.02, .12)
\\ \hline
$.1$&$10$
    &1.0(.42, 2.3)
    &.17(.05, .55)
    &.05(.02, .13)
\\ \hline
$.1$&$100$
    &6.8(1.7, 17.9)
    &1.3(.30, 5.0)
    &.05(.02, .14)
\\ \hline
\end{tabular}
\end{center}
\end{table}

The numerical experiment uses 
two full periods of $X(t)$ generated with \( N = 100 \)
and \( \Delta t = 4 \pi / N \); i.e.,
discrete time points equally spaced over the interval $[0, 4 \pi]$.
For \( k = 1 , \ldots , N \) the measurements \( z_k \) were simulated by
\(
z_k = X_2 ( t_k ) + v_k \; .
\)
In order to test the robustness of the $\ell_1$ model to measurement noise
containing outlier data,
we generate $v_k$ as a
mixture of two normals with
$p$ denoting the fraction of outlier contamination; i.e.,
\begin{equation}
\label{SimulationModel}
v_k \sim (1 - p) \B{N} ( 0, 0.25 ) + p \B{N} ( 0, \phi ) \; .
\end{equation}
This was done for $p \in \{ 0, \; 0.1 \}$
and $\phi \in \{ 1, 4, 10, 100 \}$.
The model for the mean of $z_k$ given $x_k$ is
\(
h_k ( x_k ) =  ( 0 , 1 ) x_k = x_{2,k} \; .
\)
Here \( x_{2,k} \) denotes the second component of \( x_k \).
The model for the variance of $z_k$ given $x_k$ is $R_k = 0.25$.
This simulates a lack of knowledge of
the distribution for the outliers; i.e, $p \B{N} ( 0, \phi )$.
Note that we are recovering estimates for the smooth function
$ - \sin(t) $ and its derivative $ - \cos (t) $ using noisy
measurements (with outliers) of the function values.

\begin{figure}
 \hspace{-.75cm}
 \begin{center}
{\includegraphics[scale=0.35]{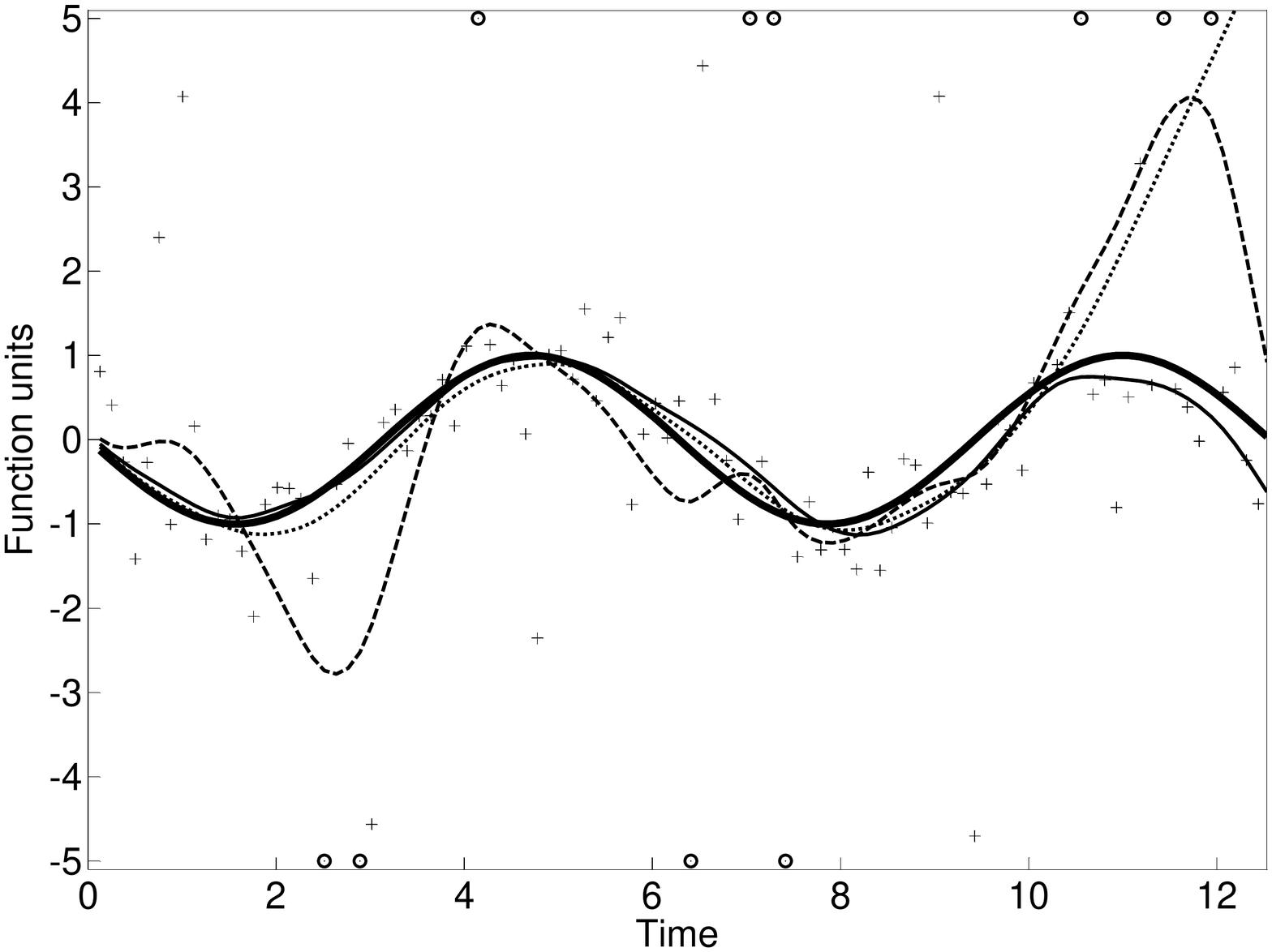}}
\end{center}
\hspace{-1cm}
\caption{
        \label{robustSpline}
        Simulation:
        measurements (+),
        outliers (o) (absolute residuals more than three standard deviations),
        true function (thick line),
        $\ell_1$-Laplace estimate (thin line),
        Gaussian estimate (dashed line),
        Gaussian outlier removal estimate (dotted line)
}
\end{figure}

We simulated 1000 realizations of the sequence \( \{ z_k \} \)
keeping the ground truth fixed,
and for each realization, and each estimation method, we computed
the corresponding state sequence estimate \( \{ \hat{x}_k \} \).
The Mean Square Error (MSE) corresponding to such an estimate is defined by
\begin{equation}
\label{MSE}
\R{MSE} =
	\frac{1}{N} \sum_{k=1}^N
		[ x_{1,k} - \hat{x}_{1,k} ]^2
		+
		[ x_{2,k} - \hat{x}_{2,k} ]^2
\; ,
\end{equation}
where \( x_k = X ( t_k ) \).
In Table~\ref{SimulationResults},
the Gaussian Kalman Filter is denoted by (GKF),
the Iterated Gaussian Smoother (IGS),
and the Iterated $\ell_1$-Laplace Smoother (ILS).
For each of these estimation techniques,
each value of $p$,
and each value of $\phi$,
the corresponding table entry is the median MSE
followed by the centralized 95\% confidence interval for the MSE.
For this problem, the model functions \( \{ g_k ( x_{k-1} ) \} \)
and \( \{ h_k ( x_k ) \} \) are linear so
the iterated smoothers IGS and ILS only require one iteration
to estimate the sequence \( \{ \hat{x}_k \} \).

Note the $\ell_1$-Laplace smoother
performs nearly as well as the Gaussian smoother
at the nominal conditions (\( p = 0 \)).
The $\ell_1$-Laplace smoother performs better and more consistently
in cases with data contamination ( \( p \geq .1 \) and \( \phi \geq 1 \) ).
It is also apparent that the smoothers perform better
than the filters.

Outlier detection and removal followed by refitting
is a simple approach to robust estimation and
can be applied to the smoothing problem.
An inherent weakness of this approach is that the
outlier detection is done using an initial fit
which assumes outliers are not present.
This can lead to good data being classified as outliers
and result in over fitting the remaining data.
An example of this is
illustrated in Figure~\ref{robustSpline}
which plots the estimation results for a realization of \( \{ z_k \} \)
where $p = 0.1$ and $\phi = 100$.
Outlier removal also makes critical review of the model more difficult.
A robust smoothing method with a consistent model,
such as the $\ell_1$-Laplace smoother,
does not suffer from these difficulties.

\subsubsection{Stochastic Nonlinear Process Example}

We now illustrate the behavior of the $\ell_1$-Laplace smoother on the Van Der Pol
Oscillator described in Section \ref{VanDerPol1}. The numerical experiment we describe
is taken from \cite[Section VI]{Aravkin2011tac}.
The corresponding nonlinear differential equation is
\[
\dot{X}_1 (t) = X_2 (t)
\hspace{.5cm} \mbox{and} \hspace{.5cm}
\dot{X}_2 (t) = \mu [ 1 - X_1 (t)^2 ] X_2 (t) - X_1 (t)
\; .
\]
Given \( X( t_{k-1} ) = x_{k-1} \) the Euler approximation for
\( X ( t_{k-1} + \Delta t ) \) is
\[
g_k ( x_{k-1} ) = \left( \begin{array}{cc}
	x_{1,k-1} + x_{2,k-1} \Delta t
	\\
	x_{2,k-1} + \{ \mu [ 1 - x_{1,k}^2 ] x_{2,k} - x_{1,k} \} \Delta t
\end{array} \right) \; .
\]
For this simulation, the `ground truth' is obtained from
a stochastic Euler approximation
of the Van der Pol oscillator.
To be specific,
with \( \mu = 2 \), \( N = 164 \) and \( \Delta t = 16 / N \),
the ground truth state vector \( x_k \) at time \( t_k = k \Delta t \)
is given by \( x_0 = ( 0 , -0.5 )^\R{T} \) and
for \( k = 1, \ldots , N \),
\begin{equation}
\label{VanDerPolTruth}
x_k = g_k ( x_{k-1} ) + w_k \; ,
\end{equation}
where $\{ w_k \}$ is a realization of
independent Gaussian noise with variance $0.01$.
Our model for state transitions (\ref{IntroGaussModel})
uses \( Q_k = 0.01 \; I \)
for \( k > 1 \), and so is identical to the model used
to simulate the ground truth \( \{ x_k \} \).
Thus, we have precise knowledge of the process that generated the
ground truth \( \{ x_k \} \).
The initial state \( x_0 \) is imprecisely specified by setting
\( g_1 ( x_0 ) = ( 0.1 , -0.4 )^\R{T} \neq x_0 \)
with corresponding variance \( Q_1 = 0.1 \; I \).
\begin{table}
\caption{
\label{SimulationResultsNL}
Median MSE over 1000 runs and
confidence intervals containing 95\% of MSE results
}
\begin{center}
\begin{tabular}{|c|c|c|c|}\hline
$\B{p}$&$\B{\phi}$
	& IGS
	& ILS
\\ \hline
$0$&$-$
    &  0.07 (0.06, 0.08)
    &  0.07 (0.06, 0.09)
\\ \hline
$.1$&$10$
    &  0.07 (0.06, 0.10)
    &  0.07 (0.06, 0.09)
    \\ \hline
$.2$&$10$
    &  0.08 (0.06, 0.11)
    &  0.08 (0.06, 0.11)
\\ \hline
$.3$&$10$
    &  0.08 (0.06, 0.11)
    &  0.08 (0.06, 0.11)
\\ \hline
$.1$&$100$
    &  0.10 (0.07, 0.14)
    &  0.07 (0.06, 0.10)
\\ \hline
$.2$&$100$
    &  0.12 (0.07, 0.40)
    &  0.08 (0.06, 0.11)
\\ \hline
$.3$&$100$
    &  0.13 (0.09, 0.64)
    &  0.08 (0.07, 0.10)
\\ \hline
$.1$&$1000$
    &  0.17 (0.11, 1.50)
    &  0.08 (0.06, 0.11)
\\ \hline
$.2$&$1000$
    &  0.21 (0.14, 2.03)
    &  0.08 (0.06, 0.11)
\\ \hline
$.3$&$1000$
    &  0.25 (0.17, 2.66)
    &  0.09 (0.07, 0.12)
\\ \hline
\end{tabular}
\end{center}
\end{table}

\begin{figure*}
{\hskip -1in \includegraphics[scale=0.5]{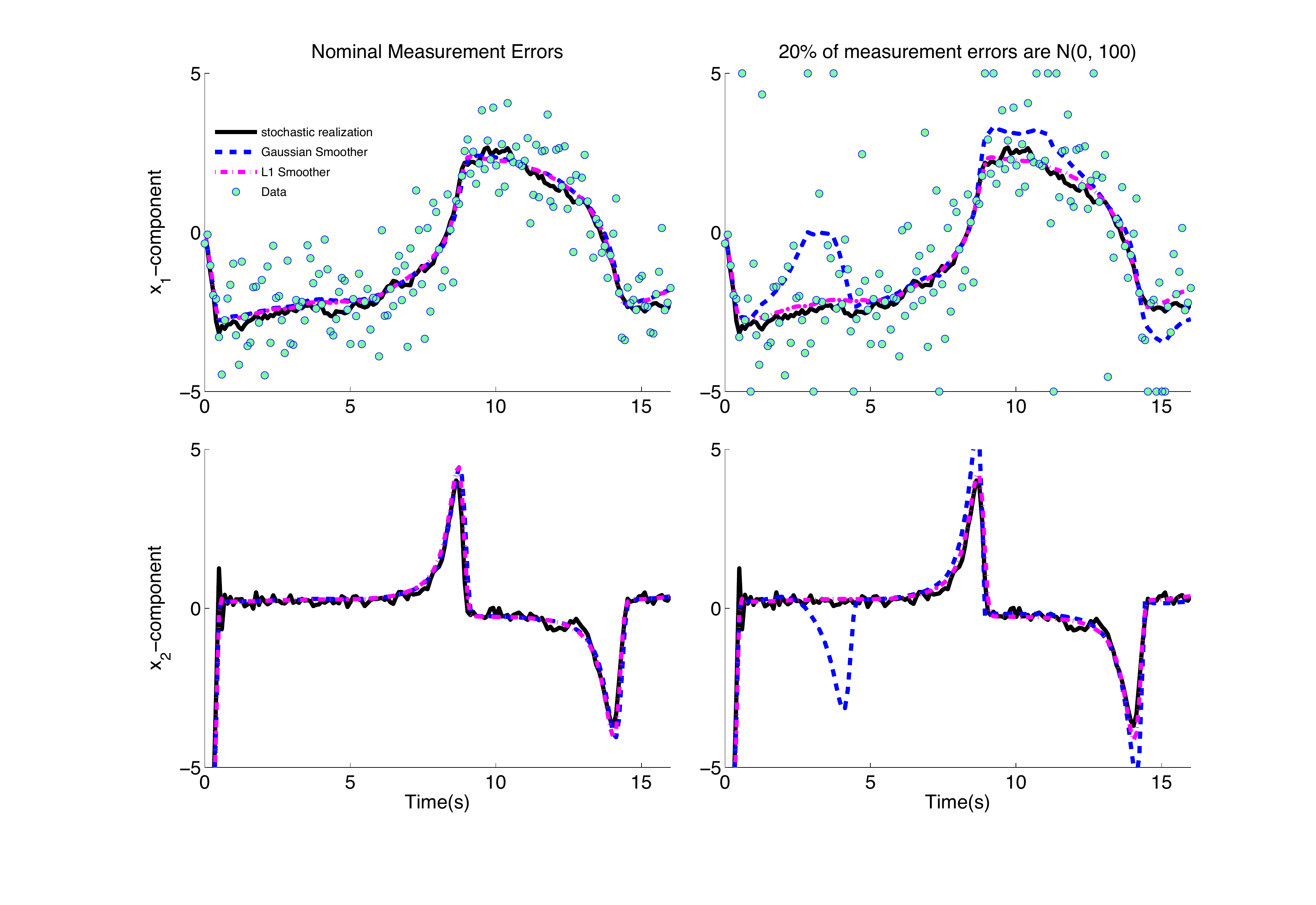}}
\vspace{ -.3in}
\caption{
	\label{nlSim}
       The left two panels show estimation of $x_{1},$ (top)
       and  $x_2$ (bottom) with errors from the nominal model.
       The stochastic realization is represented by a thick black line;
       the Gaussian smoother is the blue dashed line, and
       the $\ell_1$-smoother is the magenta dash-dotted line.
       Right two panels show the same stochastic realization but
       with measurement errors now from $(p, \phi) = (.2, 100)$.
	Outliers appear on the top and bottom boundary in the top right panel.
	}
\end{figure*}

For \( k = 1 , \ldots , N \) the measurements \( z_k \) were simulated by
\(
z_k = x_{1,k} + v_k \; .
\)
The measurement noise $v_k$ was generated as follows:
\begin{equation}
v_k \sim (1 - p) \B{N} ( 0, 1.0 ) + p \B{N} ( 0, \phi ) \; .
\end{equation}
This was done for
$p \in \{ 0, 0.1, 0.2, 0.3 \}$ and
$ \phi \in \{ 10, 100, 1000 \}$.
The model for the mean of $z_k$ given $x_k$ is
\(
h_k ( x_k ) =  ( 1 , 0 ) x_k = x_{1,k}
\).
As in the previous simulation, we simulated a lack of knowledge of
the distribution for the outliers; i.e, $p \B{N} ( 0, \phi )$.
In (\ref{IntroGaussModel}),
the model for the variance of $z_k$ given $x_k$ is $R_k = 1.0$.

We simulated 1000 realizations of the
ground truth state sequence \( \{ x_k \} \)
and the corresponding measurement sequence \( \{ z_k \} \).
For each realization, we computed
the corresponding state sequence estimate \( \{ \hat{x}_k \} \)
using both the IGS and IKS procedures.
The Mean Square Error (MSE) corresponding to such an estimate is defined by
equation~(\ref{MSE}),
where \( x_k \) is given by equation~(\ref{VanDerPolTruth}).
The results of the simulation appear in Table~\ref{SimulationResultsNL}.
As the proportion and variance of the outliers increase,
the Gaussian smoother degrades, but
the $\ell_1$-Laplace smoother is not affected.

Figure \ref{nlSim} provides a visual illustration
of one realization \( \{ x_k \} \)
and its corresponding estimates \( \{ \hat{x}_k \} \).
The left two panels demonstrate that,
when no outliers are present,
both the IGS and ILS generate accurate estimates.
Note that we only observe the first component of the state
and that the variance of the observation is relatively large
(see top two panels).
The right two panels show what can go wrong when outliers are present.
The Van der Pol oscillator can have sharp peaks
as a result of the nonlinearity in its process model,
and outliers in the measurements can `trick' the IGS into these modes
when they are not really present.
In contrast, the Iterated $\ell_1$-Laplace Smoother avoids this problem.

\subsection{Further Extensions with Log-Concave Densities}

Let us step back for a moment and examine a theme common to all of the variations on the
Kalman smoother that we have examined thus far and 
compare
the objective functions in \eqref{l1MAP2}, \eqref{fullLS}, \eqref{fullNLLS}, \eqref{fullLSc},
and \eqref{conNLLS}. 
In all cases, the objective function takes the form
\begin{equation} \label{elqp prob}
\sum_{k=1}^N V_k\left(h( x_k)-z_k;R_k\right)  +  J_k\left(x_k - g( x_{k-1});Q_k \right)\;,
\end{equation}
where the mappings $V_k$ and $J_k$ are associated with log-concave densities of the form
\[
p_{v,k}(z)\propto \exp\left(-V_k(z:R_k)\right))\quad\mbox{ and }\quad 
p_{w,k}(x)\propto\exp\left(-J_k(x;Q_k)\right)
\]
with $p_{v,k}$ and $p_{w,k}$ having covariance matrices $R_k$ and $Q_k$, respectively.
The choice of the penalty functions $V_k$ and $J_k$ reflect the underlying model for distribution
of the observations and the state, respectively. In many applications, the functions $V_k$ and $J_k$ 
are a members
of the class
of extended piecewise linear-quadratic penalty functions.

\subsubsection{Extended Linear-Quadratic Penalties}

\begin{definition}\label{def:eplq}
For a nonempty polyhedral
set $U \subset \mB{R}^m$ and a symmetric positive-semidefinite
matrix $M\in \mB{R}^{m\times m}$ (possibly $M =0$), define the function
$\theta_{U, M}: \mB{R}^m \rightarrow \{\mB{R} \cup \infty\}:=\overline{\mB{R}}$ by
\begin{equation}
\label{PLQbasic}
\theta_{U, M}(w) := \sup_{u \in U}\left\{\langle u,w \rangle -
\frac{1}{2}\langle u, Mu \rangle\right\}.
\end{equation}
Given and injective matrix $B\in\mB{R}^{m\times n}$ and a vector $b\in\mB{R}^m$, define
$\rho: \mB{R}^n \rightarrow \overline{\mB{R}}$
as $\theta_{U,M}(b + By)$:
\begin{equation}\label{EPLQbasic}
\begin{array}{rcl}
\rho_{U, M, b, B}(y) &:=& 
\sup_{u \in U}
\left\{ \langle u,b + By \rangle - \frac{1}{2}\langle u, Mu
\rangle \right\} 
\end{array}\;.
\end{equation}
All functions of the type specified in \eqref{PLQbasic} are called {\em piecewise linear-quadratic} (PLQ)
penalty functions,
and those of the form \eqref{EPLQbasic} are called {\em extended piecewise linear-quadratic}
(EPLQ) penalty functions.
\end{definition}

\begin{remark}
PLQ penalty functions are extensively studied by Rockafellar and Wets in \cite{RTRW}. 
In particular, they present a full duality theory for optimizations problems based on these functions.
\end{remark}

It is easily seen that the penalty functions arising from both the Gaussian and $\ell_1$-Laplace 
distributions come from this EPLQ class. But so do other important densities such as the Huber
and Vapnik densities.
\smallskip

\begin{figure} \label{HuberVapnikFig}
  \begin{center} 
  \begin{tabular}{cc} 
{  \includegraphics[scale=0.29]{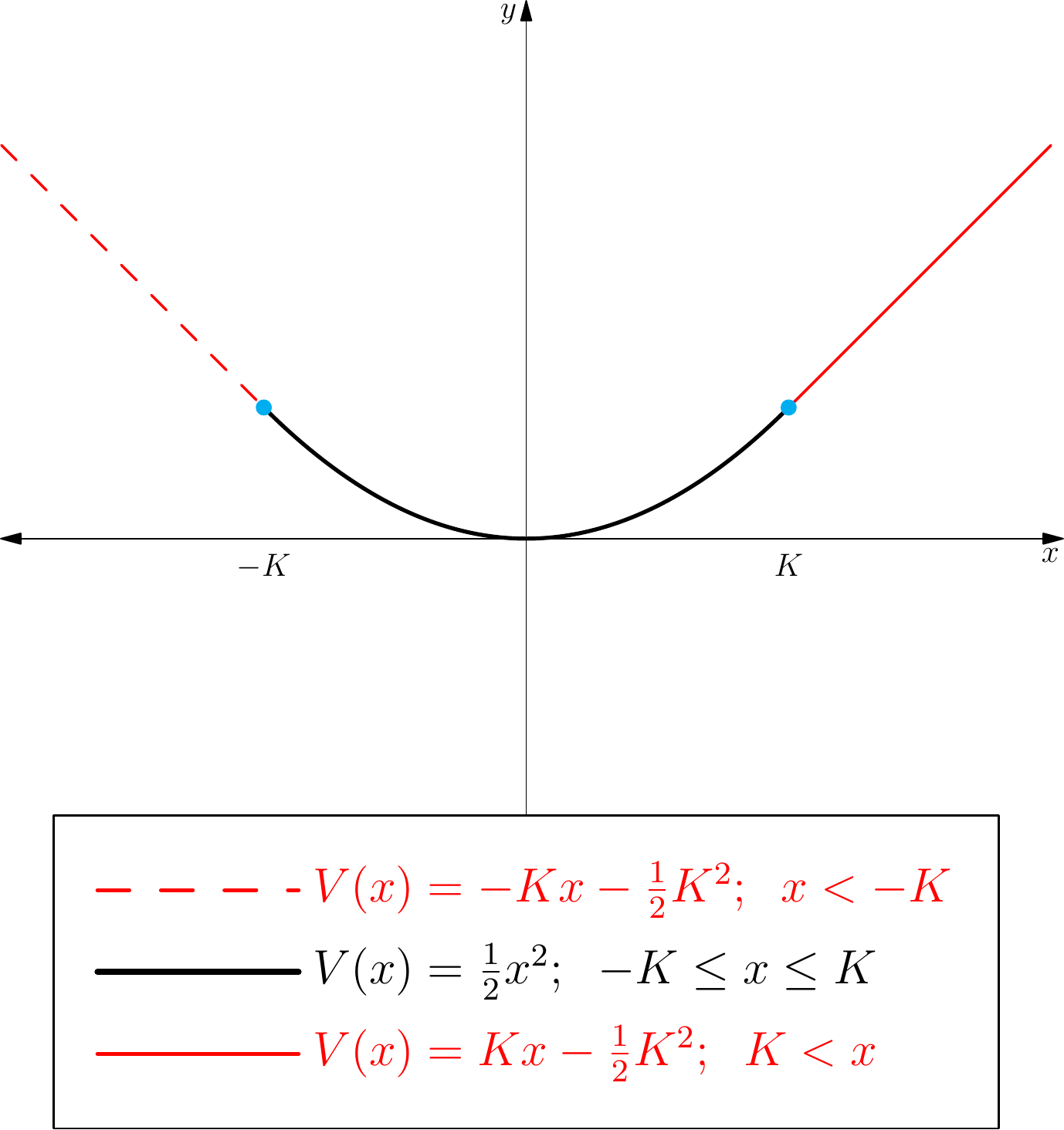}}
\hspace{.1in}
{\includegraphics[scale=0.36]{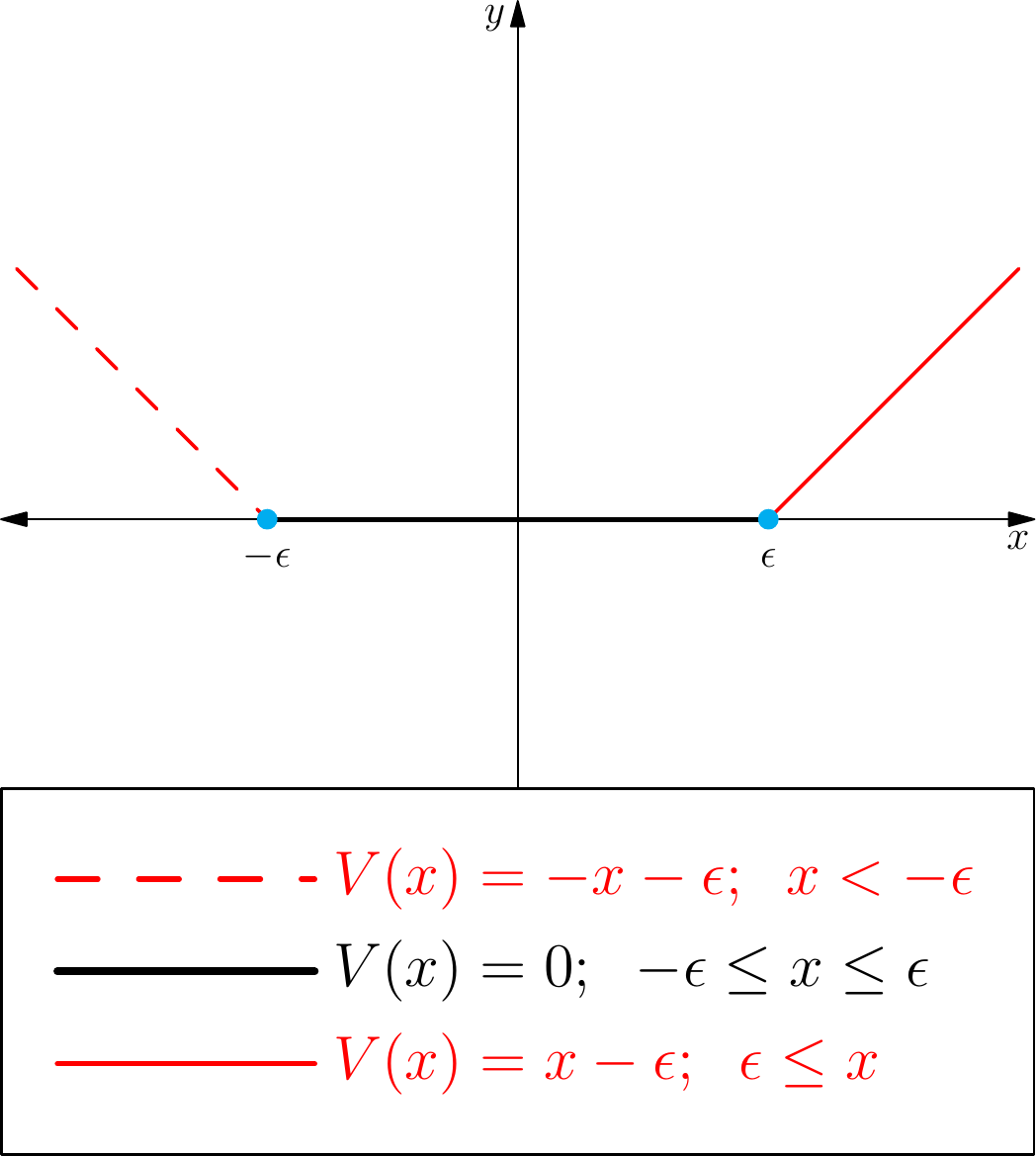}}
   \end{tabular}
    \caption{Huber (left) and Vapnik (right) Penalties}
     \end{center}
\end{figure}


\noindent
{\bf Examples:} 
The $\ell_2$, $\ell_1$, Huber, and Vapnik
penalties are representable in the notation of Definition
\ref{def:eplq}.
\begin{enumerate}
\item $L_2$: Take $U = \mB{R}$, $M = 1$, $b = 0$, and $B = 1$. We obtain
\( \displaystyle \rho(y) = \sup_{u \in \B{R}}\left\langle uy -
\frac{1}{2}u^2 \right\rangle\;. \) The function inside the $\sup$ is
maximized at $u = y$, whence $\rho(y) = \frac{1}{2}y^2$.

\item $\ell_1$: Take $U = [-1, 1]$, $M = 0$, $b = 0$, and $B = 1$. We obtain
\( \displaystyle \rho(y) = \sup_{u \in [-1, 1]}\left\langle
uy\right\rangle\;. \) The function inside the $\sup$ is maximized by
taking $u = \R{sign}(y)$, whence $\rho(y) = |y|$.

\item Huber: Take $U = [-K, K]$, $M = 1$, $b = 0$, and $B = 1$.
We obtain \( \displaystyle \rho(y) = \sup_{u \in [-K,
K]}\left\langle uy - \frac{1}{2}u^2 \right\rangle\;. \) Take the
derivative with respect to $u$ and consider the following cases:
\begin{enumerate}
\item If $ y < -K $, take $u = -K$ to obtain
$-Ky  -\frac{1}{2}K^2$.
\item If $-K \leq y \leq K$, take $u = y$ to obtain
$\frac{1}{2}y^2$.
\item If $y > K $, take $u = K$ to obtain
a contribution of $Ky -\frac{1}{2}K^2$.
\end{enumerate}
This is the Huber penalty with parameter $K$, shown in the left
panel of Fig. 1.

\item Vapnik: take $U = [0,1]\times[0,1]$,
$M = \left[\begin{smallmatrix}0 & 0\\0 & 0
\end{smallmatrix}\right]$, $B = \left[ \begin{smallmatrix} 1\\-1
\end{smallmatrix} \right]$, and $b = \left[ \begin{smallmatrix}
-\epsilon \\-\epsilon \end{smallmatrix} \right]$, for some $\epsilon
> 0$. We obtain \( \rho(y) = \sup_{u_1, u_2 \in [0,1]} \left\langle
\begin{bmatrix}
y-\epsilon\\
-y-\epsilon
\end{bmatrix},
\begin{bmatrix}
u_1\\
u_2
\end{bmatrix}
\right\rangle . \)
We can obtain an explicit representation by considering three cases:
\begin{enumerate}
\item If $|y| < \epsilon$, take $u_1 = u_2 = 0$. Then $\rho(y) = 0$.
\item If $y > \epsilon$, take $u_1 = 1$ and $u_2 = 0$. Then
$\rho(y) = y - \epsilon$.
\item If $y < -\epsilon$, take $u_1 = 0$ and $u_2 = 1$. Then
$\rho(y) = -y - \epsilon$.
\end{enumerate}
This is the Vapnik penalty with parameter $\epsilon$, shown in the
right panel of Fig.  \ref{HuberVapnikFig}.
\end{enumerate}
\smallskip

\subsubsection{PLQ Densities}

We caution that not every EPLQ function is the negative log of a density function. For an
ELQP function $\rho$ to be associated with a density, the function
$\exp(-\rho(x))$ must be integrable on $\mB{R}^n$. The integrability of $\exp(-\rho(x))$
can be established under a coercivity hypothesis.

\begin{definition}\label{defn:coercive}
A function $\rho:\mB{R}^n\rightarrow \mB{R}\cup\{+\infty\}=\eR$ is said to be coercive
(or $0$-coercive) if $\lim_{\Vert x\Vert\rightarrow \infty}\rho(x)=+\infty$.
\end{definition}

Since the functions $\rho_{U,M,b,B}$ defined in \eqref{EPLQbasic} are not necessarily
finite-valued, their calculus must be treated with care. An important tool in this regard
is the essential dominion. The
essential domain of $\rho:\mB{R}^n\rightarrow \eR$ is the set
\[
\R{dom}(\rho):=\left\{x\, :\, \rho(x)<+\infty\right\}.
\]
The affine hull of $\R{dom}(\rho)$ is the smallest affine set containing $\R{dom}(\rho)$, where
a set is affine if it is the translate of a subspace.

\begin{theorem}\cite[Theorem 6]{RobustSparseArxiv}
\label{PLQIntegrability} (PLQ Integrability). 
Let $\rho:=\rho_{U,M,b,B}$ be defined as in \eqref{EPLQbasic}.
Suppose $\rho(y)$ is
coercive, and let $n_{\R{aff}}$ denote the dimension of
$\R{aff}(\R{dom}\; \rho)$. Then the function $f(y) = \exp(-\rho(y))$
is integrable on $\R{aff}(\R{dom}\; \rho)$ with the
$n_{\R{aff}}$-dimensional Lebesgue measure.
\begin{flushright}
$\blacksquare$
\end{flushright}\end{theorem}

\begin{theorem}\cite[Theorem 7]{RobustSparseArxiv}
\label{coerciveRho} (Coercivity of $\rho$). 
The function $\rho_{U,M,b,B}$ defined in \eqref{EPLQbasic} is coercive if
and only if $[B^\R{T}\mathrm{cone}(U)]^\circ = \{0\}$.
\begin{flushright}
$\blacksquare$
\end{flushright}\end{theorem}

If  $\rho:=\rho_{U,M,b,B}$ is coercive, then, by
Theorem \ref{PLQIntegrability}, then the function $f(y) = \exp(-\rho(y))$
is integrable on $\R{aff}(\R{dom}\; \rho)$ with the
$n_{\R{aff}}$-dimensional Lebesgue measure. If we define
\begin{equation}
\label{PLQdensity} \B{p}(y) =
\begin{cases}
c_1^{-1}\exp(- \rho(y)) & y \in \R{dom}\; \rho\\
0 & \R{else},
\end{cases}
\end{equation}
where
\[
c_1 = \left(\int_{y \in \R{dom}\; \rho} \exp(-\rho(y))dy\right),
\]
and the integral is with respect to the Lebesgue measure with
dimension $n_{\R{aff}}$, 
then $\B{p}$ is a probability density on $\R{dom}( \rho)$. 
We call these PLQ densities.

\subsubsection{PLQ Densities and Kalman Smoothing}

We now show how to build up the penalty functions $V_k$ and $J_k$ in \eqref{elqp prob} 
using PLQ densities. We will do this for the linear model \eqref{IntroGaussModel}-\eqref{IntroLinearModel} 
for simplicity.
The nonlinear case can be handled as before by applying the Gauss-Newton
strategy to the underlying convex composite function.

Using the notion given in \eqref{defs}, the linear model \eqref{IntroGaussModel}-\eqref{IntroLinearModel}
can be written as  
\begin{equation}
\label{fullStat}
\begin{array}{lll}
w
&=&
Gx + \B{w}\\
z 
&=&
Hx + \B{v}\;.
\end{array}
\end{equation}

A general Kalman smoothing problem can be specified 
by assuming that the noises $\B{w}$ and $\B{v}$ in the model~\eqref{fullStat}
have PLQ densities with means $0$, variances  $Q$ and $R$
\eqref{defs}. Then, for suitable
$\{U^{w}_{k}, M^{w}_{k},b^{w}_{k},B^{w}_{k}\}$ and $\{U^{v}_{k}, M^{v}_{k},b^{v}_{k},B^{v}_{k}\}$, 
we have 
\begin{equation}
\label{kalmanDensities}
\begin{aligned}
\B{p}(w)&\propto \exp(-\theta_{U^w, M^w}(b^w + B^w Q^{-1/2}w)) \\
 \B{p}(v) &\propto \exp(-\theta_{U^v, M^v}(b^v + B^v R^{-1/2}v))\; ,
\end{aligned}
\end{equation}
where
\[
\begin{aligned}
U^w&=\prod_{k=1}^NU^{w}_{k}\subset\mB{R}^{nN}\\
U^v&=\prod_{k=1}^NU^{v}_{k}\subset\mB{R}^{M}
\end{aligned}\ ,
\quad
\begin{aligned}
M^w&=\R{diag}(\{M^{w}_{k}\})\\
M^v&=\R{diag}(\{M^{v}_{k}\})
\end{aligned}\ ,
\quad
\begin{aligned}
B^w&=\R{diag}(\{B^{w}_{k}\})\\
B^v&=\R{diag}(\{B^{v}_{k}\})\\
b^w&=\R{vec}(\{b^{w}_{k}\})\\
b^v&=\R{vec}(\{b^{v}_{k}\})
\end{aligned}.
\]
Then the MAP estimator
for $x$ in the model~\eqref{fullStat} is 
\begin{equation}
\arg \min_{x\in \mB{R}^{nN}} 
\left\{
\begin{aligned}
\label{PLQsubproblem} 
&\theta_{U^w,M^w}(b^w + B^w Q^{-1/2}(Gx - w)) \\
&+ \theta_{U^v, M^v}(b^v + B^vR^{-1/2}(Hx - z))
\end{aligned}
\right\}\;. 
\end{equation} 

Note that since $w_k$ and $v_k$ are independent,
problem~\eqref{PLQsubproblem} is decomposable into a sum of 
terms analogous to~\eqref{elqp prob}.
This special structure follows from the 
block diagonal structure of $H, Q, R, B^v, B^w$, 
the bidiagonal structure of $G$, and 
the product structure of sets $U^w$ and $U^v$, 
and is key in proving the linear complexity of the
solution method we propose.

\subsubsection{Solving the Kalman Smoother Problem with PLQ Densities}

Recall that, when the sets $U^w$ and $U^v$ are polyhedral, 
(\ref{PLQsubproblem}) is an Extended
Linear Quadratic program (ELQP), described in \cite[Example
11.43]{RTRW}. 
We solve~\eqref{PLQsubproblem} by working directly  
with its associated Karush-Kuhn-Tucker (KKT) system. 
\begin{lemma}\cite[Lemma 3.1]{RobustSparseArxiv}
\label{lem:KKT}
Suppose that the sets $U^w_k$ and $U^v_k$
are polyhedral, that is, they can be given the representation
\[
U^w_k = \{u|(A^w_k)^Tu \leq a^w_k \}, \quad U^v_k = \{u|(A^v_k)^Tu\leq a^v_k\}\;.
\]
Then the first-order necessary and sufficient conditions 
for optimality in~\eqref{PLQsubproblem} are given by 
\begin{equation}
\label{PLQFinalConditions}
\begin{array}{lll}
&\begin{array}{llllll}
0 &=& (A^w)^\R{T}u^w + s^w - a^w\;;&&  0=  (A^v)^\R{T}u^v + s^v - a^v\\
0 &=& (s^w)^\R{T}q^w\;; && 0= (s^v)^\R{T}q^v
\end{array}\\\\
&\begin{array}{llllll}
0 &=& \tilde b^w + B^w Q^{-1/2}G{x} -  M^w{u}^w - A^wq^w\\
0 &=& \tilde b^v - B^v R^{-1/2}H{x} - M^v{u}^v -  A^v q^v\\
0 &=& G^\R{T}Q^{-\R{T}/2}(B^w)^\R{T} u^w -
H^\R{T}R^{-\R{T}/2}(B^v)^\R{T} u^v\\
0 &\leq& s^w, s^v, q^w, q^v. 
\end{array}
\end{array}\;,
\end{equation}
where $\tilde b^w = b^w - B^wQ^{-1/2}w$ and 
$\tilde b^v = b^v - B^vR^{-1/2}z$.
\begin{flushright}
$\blacksquare$
\end{flushright}\end{lemma}

We propose solving the KKT conditions~\eqref{PLQFinalConditions}
by an Interior Point (IP) method.  
IP methods work by applying a damped Newton iteration to a
relaxed version of (\ref{PLQFinalConditions}) where the relaxation is to the 
complementarity conditions. Specifically, we replace the complementarity conditions by
\[
\begin{array}{lll}
(s^w)^\R{T}q^w = 0 & \rightarrow & Q^wS^w\B{1} - \mu\B{1} = 0 \\
(s^v)^\R{T}q^v = 0 & \rightarrow & Q^vS^v\B{1} - \mu\B{1} = 0\;,
\end{array}
\]
where $Q^w, S^w, Q^v, S^v$ are diagonal matrices
with diagonals $q^w, s^w, q^v, s^v$ respectively. 
The parameter $\mu$ is aggressively decreased to $0$ as the IP 
iterations proceed. Typically, no more than 10 or 20 iterations 
of the relaxed system are required to obtain a solution of~\eqref{PLQFinalConditions},
and hence an optimal solution to~\eqref{PLQsubproblem}.
The following theorem shows that the computational effort required (per IP iteration)
is linear in the number of time steps whatever PLQ density 
enters the state space model. 

\begin{theorem}\cite[Theorem 3.2]{RobustSparseArxiv}
\label{thm:PLQsmoother}
(PLQ Kalman Smoother Theorem) Suppose that all $w_k$ and $v_k$ in
the Kalman smoothing model~\eqref{IntroGaussModel}-\eqref{IntroLinearModel} come from PLQ
densities that satisfy 
$\mathrm{Null}(M)\cap U^{\infty}  = \{0\}$.
Then an IP method can be applied to solve~\eqref{PLQsubproblem}
with a per iteration computational complexity of $O(Nn^3 + Nm)$.
\begin{flushright}
$\blacksquare$
\end{flushright}\end{theorem}

The proof, which can be found in \cite{RobustSparseArxiv}, shows that IP methods for
solving~\eqref{PLQsubproblem} 
preserve the key block tridiagonal structure of the 
standard smoother. General smoothing
estimates can therefore be computed in $O(Nn^3)$ time, 
as long as the number of IP iterations is fixed 
(as it usually is in practice, to $10$ or $20$).\\
It is important to observe that the motivating examples 
all satisfy the conditions 
of Theorem \ref{thm:PLQsmoother}.

\begin{corollary}\cite[Corollary 3.3]{RobustSparseArxiv}
\label{cor:examples}
The densities corresponding to $L^1, L^2$, Huber, 
and Vapnik penalties 
all satisfy the hypotheses of Theorem \ref{thm:PLQsmoother}. 
\end{corollary}
{\bf Proof:}
We verify that $\mathrm{Null}(M) \cap \mathrm{Null}(A^\R{T}) = 0$
 for each of the four penalties. 
 In the $L^2$ case, $M$ has full rank.
 For the $L^1$, Huber, and Vapnik penalties, the respective
sets $U$ are bounded, so $U^{\infty}= \{0\}$.

\subsubsection{Numerical example: Vapnik penalty and functional recovery}

\begin{figure*}
\begin{center}
{\includegraphics[scale=0.6]{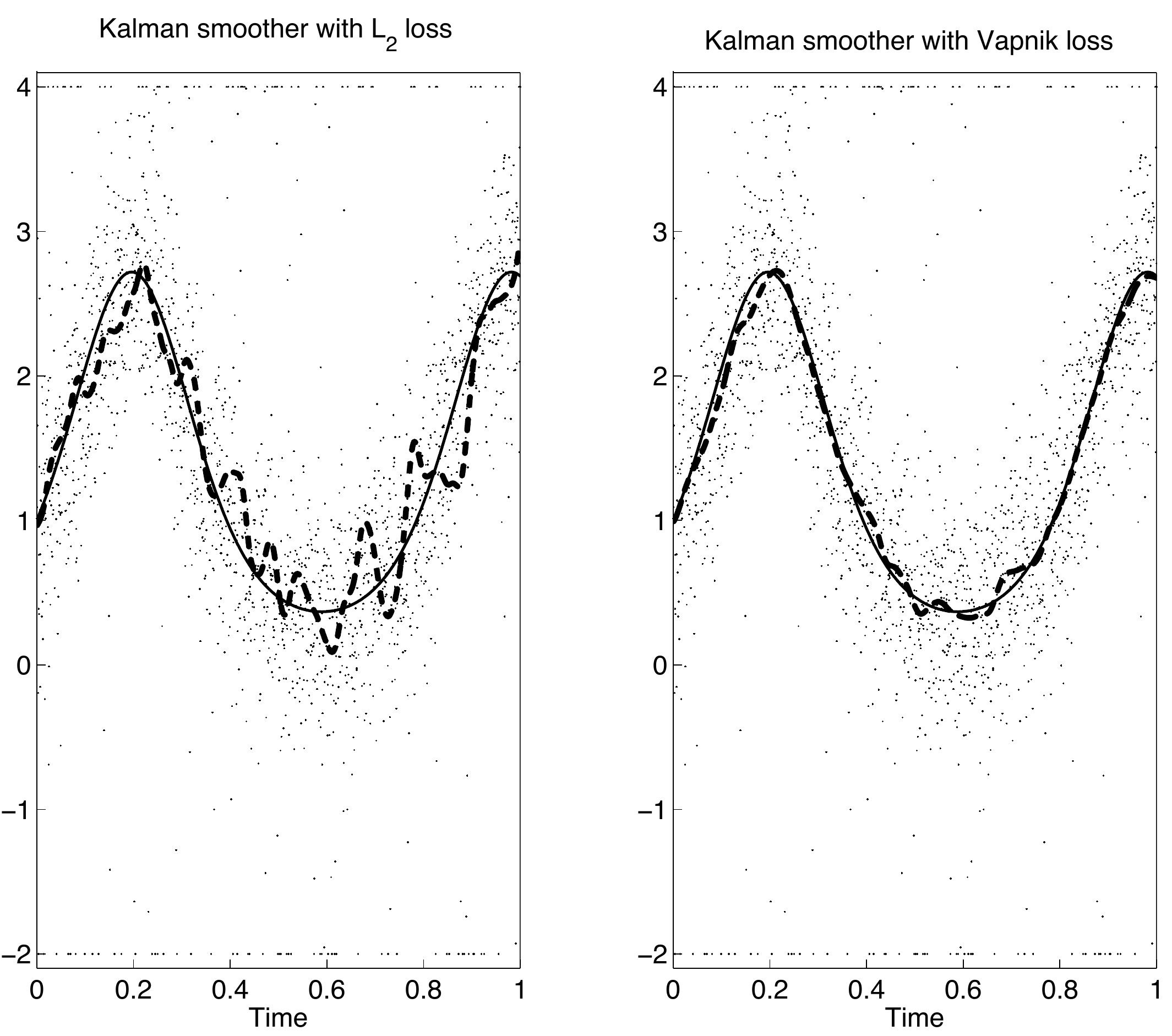}}
\end{center}
\caption{
    \label{FigLtwoVap}
    Simulation:
    measurements ($\cdot$) with outliers
    plotted on axis limits ($4$ and $-2$),
    true function (continuous line),
    smoothed estimate using either the quadratic loss (dashed line, left panel)
    or the Vapnik's $\epsilon$-insensitive loss (dashed line, right panel)
}
\end{figure*}

In this section we present a numerical example to
illustrate the use of the Vapnik penalty (see Figure~\ref{HuberVapnikFig}) in the Kalman smoothing context, 
for a functional recovery application.

We consider the following function
\begin{equation*}
f(t)=\exp\left[\sin(8t)\right]
\end{equation*}
taken from \cite{Dinuzzo07}. Our aim is to reconstruct $f$
starting from
2000 noisy samples collected uniformly over the unit interval.
The measurement noise $v_k$ was generated using a mixture of two normals with
$p=0.1$ denoting the fraction from each normal; i.e.,
\[
v_k \sim (1 - p ) \B{N} (0, 0.25 ) + p \B{N} ( 0 , 25 ),
\]
where $\B{N}$ refers to the Normal distribution. 
Data are displayed as dots in Fig.~\ref{FigLtwoVap}.
Note that the purpose of the second component of the normal mixture
is to simulate outliers in the output data and that all the measurements exceeding
vertical axis limits are plotted on upper and lower axis limits (4 and -2) to improve readability.

\begin{figure}
\begin{center}
{\includegraphics[scale=0.3]{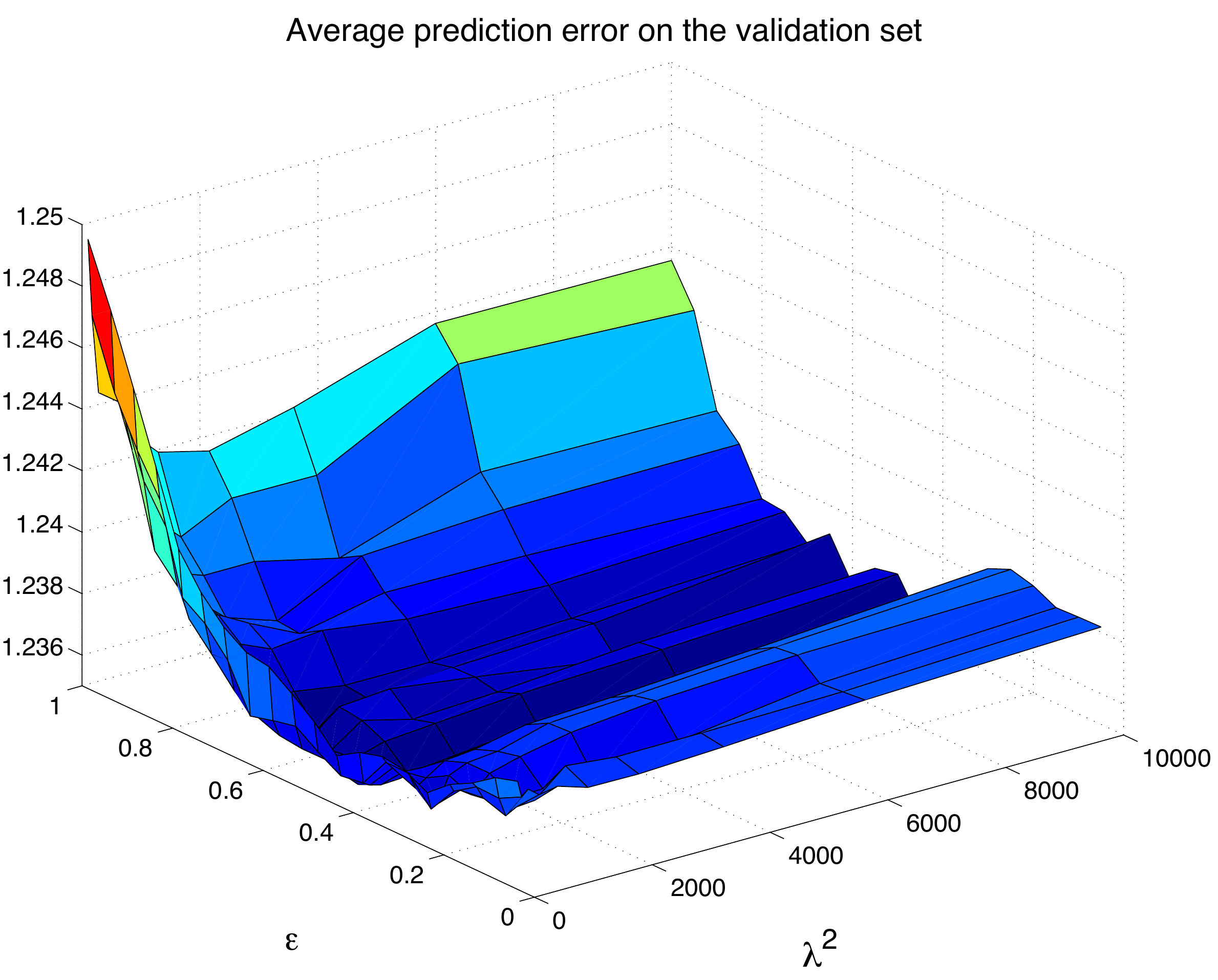}}
\end{center}
\caption{
    \label{FigVal}
    Estimation of the smoothing filter parameters using the Vapnik loss. 
    Average prediction error on the validation data set
    as a function of the variance process $\lambda^2$
    and $\epsilon$.
}
\end{figure}

The initial condition $f(0) = 1$ is assumed to be known, 
while the difference of the unknown function from the initial condition 
(i.e. $f(\cdot) - 1$) is modeled as a Gaussian process
given by an integrated Wiener process. This model captures the Bayesian
interpretation of cubic smoothing splines \cite{Wahba1990}, 
and admits a 2-dimensional state space representation
where the first component of $x(t)$, which models $f(\cdot) - 1$,
corresponds to the integral of the second state component, modelled as Brownian motion.
To be more specific, letting $\Delta t = 1/2000$, the sampled version
of the state space model (see \cite{Jaz,Oks} for details) is defined by
\begin{eqnarray*}
G_k =
\begin{bmatrix}
    1        & 0
    \\
    \Delta t & 1
\end{bmatrix}, \qquad k=2,3,\ldots,2000\\
H_k =
\begin{bmatrix}
    0        & 1
\end{bmatrix}, \qquad k=1,2,\ldots,2000
\end{eqnarray*}
with the autocovariance of $w_k$ given by
\[
Q_k = \lambda^2
\begin{bmatrix}
    \Delta t   & \frac{\Delta t^2}{2}
    \\
    \frac{\Delta t^2}{2} & \frac{\Delta t^3}{3}
\end{bmatrix}, \qquad k=1,2,\ldots,2000
\; ,
\]
where $\lambda^2$ is an unknown scale factor
to be estimated from the data.\\
The performance of two different Kalman smoothers are compared.
The first (classical) estimator uses a quadratic loss function
to describe the negative log of the measurement noise density and
contains only $\lambda^2$ as unknown parameter.
The second estimator is a Vapnik smoother relying
on the $\epsilon$-insensitive loss,
and so depends on two unknown parameters $\lambda^2$ and $\epsilon$.
In both of the cases, the unknown parameters are estimated
by means of a cross validation strategy where the 2000 measurements
are randomly split into a training and a validation set of 1300
and 700 data points, respectively. The Vapnik smoother
was implemented by exploiting the efficient computational
strategy described in the previous section, see~\cite{AravkinIFAC} for 
specific implementation details.
In this way, for each value of $\lambda^2$ and $\epsilon$
contained in a $10 \times 20$ grid on $[0.01,10000] \times [0,1]$,
with $\lambda^2$ logarithmically spaced,
the function estimate was rapidly obtained by the new smoother applied to the training set.
Then,  the relative average prediction error on the validation set was computed,
see Fig.~\ref{FigVal}.
The parameters leading to the best prediction were 
$\lambda^2=2.15\times 10^3$ and $\epsilon=0.45$, 
which give a sparse solution
defined by fewer than 400 support vectors. 
The value of $\lambda^2$
for the classical Kalman smoother was then estimated following the
same strategy described above. 
In contrast to the Vapnik penalty,
the quadratic loss does not induce any sparsity, so that, 
in this case, the number of support vectors
equals the size of the training set.\\
The left and right panels of Fig.~\ref{FigLtwoVap} display the function estimate
obtained using the quadratic and the Vapnik losses, respectively.
It is clear that the Gaussian estimate is heavily affected by the outliers. 
In contrast, as expected, the estimate coming from the Vapnik based smoother 
performs well over the entire time period, 
and is virtually unaffected by the presence of large outliers.

\section{Sparse Kalman smoothing}
\label{Sparse}

In recent years, sparsity promoting formulations and algorithms have made a tremendous
impact in signal processing, reconstruction algorithms, statistics, and inverse problems
(see e.g.~\cite{Bruckstein:2009} and the references therein). In some contexts, 
rigorous mathematical theory is available that can guarantee recovery from under-sampled
sparse signals~\cite{Donoho2006}. In addition, for many inverse problems, sparsity promoting
optimization provides a way to exploit prior knowledge of the signal class as a way to improve
the solution to an ill-posed problem, but conditions for recoverability have not yet been 
derived~\cite{Mansour11TRssma}. 

In the context of dynamic models, several sparse Kalman filters have 
been recently proposed~\cite{report:08,cskf,cskf:08, Angelosante2009}. 
In the applications considered, 
in addition to process and measurement models, the state space is also 
known to be sparse. The aim is to improve recovery by incorporating 
sparse optimization techniques. Reference~\cite{Angelosante2009} 
is very close to the work presented in this section, since they formulate
a sparsity promoting optimization problem over the whole measurement sequence
and solve it with an optimization technique shown to preserve computational efficiency.

In this section, we formulate the sparse Kalman smoothing problem as an optimization problem 
over the entire state space sequence, and suggest two new approaches for the solution of such problems. 
The first approach is based on the interior point methodology, 
and is a natural extension of the mathematics presented in earlier sections. 

The second approach is geared towards problems where the dimension $n$ 
(state at a single time point) is large. 
For this case, we propose a matrix free approach, using a different (constrained) 
Kalman smoothing formulation, together with the projected gradient method. 
In both methods, the structure of the Kalman smoothing problem is exploited to achieve computational efficiency. 

We present theoretical development for the two approaches, leaving applications and numerical 
results to future work. 

\subsection{Penalized Formulation and Interior Point Approach}

We consider only the linear smoother~\eqref{fullLS}. A straight forward way to impose sparsity 
on the state is to augment this formulation with a $1$-norm penalty: 

\begin{equation}\label{fullLSsparse}
\min_{x} f(x) :=  \frac{1}{2}\|Hx - z\|_{R^{-1}}^2 + \frac{1}{2}\|Gx - w\|_{Q^{-1}}^2 + \lambda\|Wx\|_1\;,
\end{equation}
where $W$ is a diagonal weighting matrix included for modeling convenience. For example, 
the elements of $W$ can be set to $0$ to exclude certain parts of the state dimension from 
the sparse penalty. A straightforward constrained reformulation of~\eqref{fullLSsparse} is
\begin{equation}\label{fullLSsparseCon}
\begin{aligned}
\min_{x} & \quad  \frac{1}{2}\|Hx - z\|_{R^{-1}}^2 + \frac{1}{2}\|Gx - w\|_{Q^{-1}}^2 + \lambda\B{1}^Ty\\
& \text{s.t.} \quad -y \leq Wx \leq y\;.
\end{aligned}
\end{equation}
Note that this is different from the constrained problem~\eqref{fullLSc},
because we have introduced a new variable $y$, with constraints in $x$ and $y$.
Nonetheless, an interior point approach may still be used to solve the resulting problem. 
We rewrite the constraint in~\eqref{fullLSsparseCon} using non-negative slack variables $s, r$:
\begin{equation}
\label{SparseSlack}
\begin{aligned}
Wx - y + s &= 0\\
-Wx - y + r & = 0\;,
\end{aligned}
\end{equation}
and form the Lagrangian for the corresponding system: 
\begin{equation}
\label{LagrangianSparse}
L(s, r, q, p, y, x) = x^TCx + c^Tx + \lambda\B{1}^Ty + q^T(Wx - y + s) + p^T(-Wx -y + r)\;,
\end{equation}
with $C$ as in~\eqref{hessianApprox} and $c$ as in~\eqref{defGrad}., and where $q$ and $p$ are the dual variables corresponding to the inequality constraints
$Wx \leq y$ and $-Wx \leq -y$, respectively. 
The (relaxed) KKT system is therefore given by 
\begin{equation}
\label{KKTsys}
F_\mu(s, r, q, p, y, x) :=
\left(
\begin{aligned}
 &   s - y + Wx\\ 
 &  r -  y-Wx  \\
 &  D(s)D(q)\B{1} - \mu\B{1}\\
 &  D(r)D(p)\B{1} - \mu\B{1}\\
 &  \lambda\B{1} - q - p \\
 &  Wq - Wp+ Cx + c
\end{aligned}
\right) = 0\;.
\end{equation}
The derivative matrix $F_\mu^{(1)}$ is given by
\begin{equation}
\label{NewtonSparse}
F_\mu^{(1)} 
=
\begin{bmatrix}
I & 0 & 0 & 0 & - I & W\\
0 & I & 0 & 0 & - I & - W\\
D(q) & 0 & D(s) & 0 & 0 & 0\\
0 & D(p) & 0 & D(r) & 0 & 0 & \\
0 & 0 & - I & - I & 0 & 0\\
0 & 0 & W & -W & 0 & C
\end{bmatrix}\;,
\end{equation}
and it is row equivalent to the system

\[
\begin{bmatrix}
I & 0 & 0 & 0 & - I & W\\
0 & I & 0 & 0 & - I & - W\\
0 & 0 & D(s) & 0 & D(q) & -D(q)W\\
0 & 0 & 0 & D(r) & D(p) & D(p)W  \\
0 & 0 & 0 & 0 & \Phi& -\Psi W\\
0 & 0 & 0 & 0 & 0& C +W\Phi^{-1}\left(\Phi^2 - \Psi^2\right)W
\end{bmatrix}
\]
where 
\begin{equation}
\label{PhiPsi}
\begin{aligned}
\Phi = D(s)^{-1}D(q) + D(r)^{-1}D(p)\\
\Psi = D(s)^{-1}D(q) - D(r)^{-1}D(p) \;,
\end{aligned}
\end{equation}
and the matrix $\Phi^2 - \Psi^2$ is diagonal, with 
the $ii$th entry given by $4q_ir_i$.
Therefore, the modified system preserves the structure of $C$; specifically it is symmetric, block tridiagonal, and positive definite. 
The Newton iterations required by the interior point method can therefore be carried out, with each iteration having complexity
$O(n^3N)$. 

\subsection{Constrained Formulation and Projected Gradient Approach}

Consider again the linear smoother~\eqref{fullLS}, but now impose a $1$-norm constraint rather than 
a penalty: 

\begin{equation}\label{fullLSsparseCon}
\begin{aligned}
\min_{x} f(x) :=  & \frac{1}{2}\|Hx - z\|_{R^{-1}}^2 + \frac{1}{2}\|Gx - w\|_{Q^{-1}}^2\\
\text{s.t.} \quad & \|Wx\|_1 \leq \tau\;.
\end{aligned}
\end{equation}

This problem, which equivalent to~\eqref{fullLSsparse} for certain values of $\lambda$ and $\tau$,
is precisely the LASSO problem~\cite{tibshirani:96}, and can be written 
\begin{equation}\label{LSsparseCon}
\min \frac{1}{2}x^T C x + c^Tx \quad \text{s.t.}\quad  \|Wx\|_1 \leq \tau\;.
\end{equation}
with $C\in \mathbb{R}^{nN \times nN}$ as in~\eqref{hessianApprox} and $c\in\mathbb{R}^{nN}$ as in~\eqref{defGrad}.
When $n$ is large, the interior point method proposed in the previous section may not be feasible, since it requires
exact solutions of the system 
\[
(C +W\Phi^{-1}\left(\Phi^2 - \Psi^2\right)W)x = r\;,
\] 
and the block-tridiagonal algorithm~\ref{blockAlg} requires the inversion of $n\times n$ systems. 

The problem~\eqref{LSsparseCon} can be solved without inverting such systems, using 
the spectral projected gradient method, see e.g.~\cite[Algorithm 1]{BergFriedlander:2008}. 
Specifically, the gradient $Cx + c$ must be repeatedly computed, and then 
$x^\nu - (Cx^\nu + c)$ is projected onto the set $\|Wx\|_1 \leq \tau$. 
(the word `spectral' refers to the fact that the Barzilai-Borwein line search is used to get the step length). 

In the case of the Kalman smoother, the gradient $Cx + c$ can be computed in $O(n^2N)$ time, 
because of the special structure of $C$. Thus for large systems, 
the projected gradient method that exploits the structure of $C$ affords significant savings per iteration relative to 
the interior point approach, $O(n^2N)$ vs. $O(n^3N)$, and relative to a method agnostic to the structure of $C$,  
$O(n^2N)$ vs. $O(n^2N^2)$. The projection onto the feasible set $\|Wx\|_1 \leq \tau$ can be done in $O(nN\log (nN))$ time.

\section{Conclusions}
\label{Conclusions}

In this chapter, we have presented an optimization approach to Kalman smoothing, together with a survey of 
applications and extensions. In Section~\ref{equivalence}, we showed that the recursive Kalman
filtering and smoothing algorithm is equivalent to algorithm~\ref{blockAlg}, an efficient 
method to solve block tridiagonal positive definite systems. In the following sections,
we used this algorithm as a subroutine, allowing us to present new ideas on a high level,
without needing to explicitly write down modified Kalman filtering and smoothing equations. 

We have presented extensions to nonlinear process and measurement models in Section~\ref{Nonlinear},
described constrained Kalman smoothing (both the linear and nonlinear cases) in Section~\ref{Constrained}, 
and presented an entire class of robust Kalman smoothers (derived by considering 
log-linear-quadratic densities) in Section~\ref{Robust}. For all of these applications, 
nonlinearity in the process, measurements, and constraints can be handled by a generalized
Gauss-Newton method that exploits the convex composite structure discussed in Sections~\ref{sec:ConvexComposite}
and~\ref{sec:ConvexCompositeCon}. The GN subproblem can be solved either in closed form 
or via an interior point approach; in both cases algorithm~\ref{blockAlg} was used. 
For all of these extensions, numerical illustrations have also been presented, and most
are available for public release through the \verb{ckbs{ package~\cite{ckbs}.

In the case of the robust smoothers, it is possible to extend the density 
modeling approach by considering densities outside the log-concave class~\cite{StArxiv},
but we do not discuss this work here. 

We ended the survey of extensions by considering 
two novel approaches to Kalman smoothing of sparse systems, for applications where 
modeling the sparsity of the state space sequence improves recovery. The first method
built on the readers' familiarity with the interior point approach as a tool for the constrained
extension in Section~\ref{Constrained}. The second method is suitable for large
systems, where exact solution of the linear systems is not possible. Numerical illustrations
of the methods have been left to future work.

\bibliographystyle{plain}
\bibliography{filter}

\end{document}